\documentclass[12pt]{article}
\usepackage{amsmath,amssymb,graphicx}

\newtheorem{tm}{Theorem}[section]
\newtheorem{lm}[tm]{Lemma}

\newtheorem{re}[tm]{Remark}

\newtheorem{pr}[tm]{Proposition}

\newtheorem{ex}[tm]{Example}
\newcommand{\qed}{~~\hbox{\rule{4pt}{8pt}}}

\makeatletter
 
 \@addtoreset{equation}{section}
\makeatother
\begin{document}
\setlength{\baselineskip}
{15.5pt}
\title{
An explicit effect of non-symmetry 
of random walks on the triangular lattice}
\author{
{\Large Satoshi Ishiwata
\footnote{Department of Mathematical Sciences, Faculty of Science, Yamagata University,
1-4-12, Kojirakawa, Yamagata 990-8560, Japan
(e-mail: {\tt ishiwata@sci.kj.yamagata-u.ac.jp})}
\footnote{Partially supported by Grant-in-Aid for Young Scientists (B)(No. 21740034), JSPS}, 
Hiroshi Kawabi 
\footnote{Department of Mathematics, Faculty of Science, Okayama University,
3-1-1, Tsushima-Naka, Kita-ku, Okayama 700-8530, Japan (e-mail: {\tt{kawabi@math.okayama-u.ac.jp}})}
\footnote{Partially supported by Grant-in-Aid for Young Scientists (B)(No. 23740107), JSPS}
~and Tsubasa Teruya
\footnote{The Okinawa Kaiho Bank, Ltd., 2-9-12, Kumoji, Naha, Okinawa 900-8686, Japan}
}}
\date{\today}
%
\maketitle 
%
%
\begin{abstract}
In the present paper, we study an explicit effect of non-symmetry on 
asymptotics
of the $n$-step transition probability as $n\rightarrow \infty$ for 
a class of non-symmetric random walks on the triangular lattice.
Realizing the triangular lattice into $\mathbb{R}^2$ appropriately, 
we observe that the Euclidean distance 
in $\mathbb{R}^2$ naturally appears in the asymptotics. We characterize 
this realization from a geometric view point of Kotani-Sunada's standard 
realization of crystal lattices.
As a corollary of the main theorem, we prove that 
the transition semigroup generated by 
the non-symmetric random walk
approximates the heat semigroup
generated by the usual Brownian motion on $\mathbb{R}^2$.
\end{abstract}
%
\section{Introduction}
Let $G=(V,E)$ be a locally finite, connected, oriented graph. 
Here $V$ is the set of vertices and $E$ is the set of oriented edges.
For an oriented edge $e \in E$, the {\it origin} and the {\it terminus} of $e$ are denoted by 
$o(e)$ and $t(e)$, respectively. The {\it inverse edge} of $e$ is denoted by $\overline{e}$.
A {\it random walk} on $G$ is given by a non-negative valued function $p$ on $E$ satisfying 
\begin{equation*}
\sum_{e \in E_x} p(e)=1 \quad \mbox{for all $x \in V$},
\end{equation*}
where $E_x=\left\{ e \in E \vert~ o(e)=x \right\}$. 
Here $p(e)$ is the probability that a particle at $o(e)$ moves to 
$t(e)$ along the edge $e$ in one unit time.
Then the transition probability $p(n,x,y)$ that a particle starting at $x \in V$ 
reaches $y \in V$ at time $n$ is given by
\begin{equation*}
p(n,x,y)=\sum_{c=(e_1,e_2, \ldots, e_n)} p(e_1)p(e_2) \cdots p(e_n),
\end{equation*}
where the sum is taken over all paths $c=(e_1, e_2 , \ldots , e_n)$ with
$t(e_i)=o(e_{i+1})$, $i=1, \ldots ,n-1$ and  
$o(e_1)=x$, $t(e_n)=y$.
If there exists a positive valued 
function $m_{V}$ on $V$ such that 
\begin{equation}
p(e)m_{V}(o(e))=p({\overline e})m_{V}(t(e)), \quad e\in E,
\nonumber
\label{symmetric}
\end{equation}
the random walk is said to be ($m_{V}$-)symmetric. 

A principal
theme for random walks is to investigate the properties
of $p(n,x,y)$ as $n \to \infty$. One of the most classical problems is the
recurrence-transience problem which is related to the divergence-convergence
of $\sum_{n=1}^{\infty} p(n,x,x)$, and the local central limit theorem gives 
us a useful criterion for the divergence-convergence of the series. For this reason,
this theme has been
discussed intensively in various settings by many authors. 
See Spitzer
\cite{Spitzer}, Lawler and Limic \cite{LL} and Woess \cite{Woess} for an
overview of random walks.

In particular, 
Kotani, Shirai and Sunada 
investigated
long time asymptotics of
$p(n,x,y)$
of ($m_{V}$-)symmetric random walks on a {\it crystal lattice}, a covering graph of a finite graph 
whose covering transformation group is abelian. In \cite{KS00}, as the precision of
the local central limit theorem (cf. \cite{KSS}), they established the 
asymptotic expansion 
\begin{equation}
p(n,x,y)m_{V}(y)^{-1}
\sim
a_{0}n^{-r/2} \exp \big( -\frac{d(x,y)^{2}}{4n} \big) 
\cdot
\Big (1+a_{1}(x,y)n^{-1}
+a_{2}(x,y)n^{-2}+\cdots\Big)
\label{AE}
\end{equation}
as $n \to \infty$, where $d(x,y)$ is a Euclidean distance appeared
through the {\it standard realization} of 
the graph 
into $\mathbb{R}^r$.
In their proof, spectral theoretic arguments 
due to the periodicity of the graph and 
the symmetry of the random walk play crucial roles.

Later in \cite{Ucchi, Uchiyama}, 
Uchiyama and his coauthor also
established the formula (\ref{AE}) for non-symmetric
random walks on periodic graphs (i.e., crystal lattices) in the Euclidean space by 
a probabilistic approach. 
Their result implies that the effect of the non-symmetry on the coefficient 
$a_{1}(x,y)$ highly depends not only on the underlying periodic
graph but also on the choice of the $1$-step transition probability $p$
even if the {zero mean condition} (see condition {\bf (P2)} in Section 2) 
is imposed.

In view of these results, it is a meaningful problem to determine
an explicit effect of the non-symmetry on the coefficient
$a_{1}(x,y)$ on a specific graph. In the present paper, 
we give an answer of this problem in the case of
the triangular lattice. As we will mention later, there exist non-symmetric 
random walks on the triangular lattice satisfying the zero mean condition {\bf (P2)}.
On the other hand, for example, on the square lattice and the hexagonal lattice,
we note that the zero mean condition is equivalent to
the symmetry of the random walk since these graphs are maximal abelian coverings
(see Kotani and Sunada \cite[page 842]{KS06}).
In the proof of the main theorem (Theorem \ref{Main}),
we make use of the probabilistic approach as in \cite{Ucchi, Uchiyama}
with the idea of the standard realization by Kotani and Sunada 
\cite{KS00, KS00-2}.

As a corollary of Theorem \ref{Main}, we establish 
the functional analytic central limit theorem
(Theorem \ref{FTCLT}) that 
the transition semigroup 
generated by the non-symmetric random walk on the triangular lattice 
approximates the heat semigroup
generated by the usual Brownian motion on $\mathbb{R}^2$.
It says that the effect of the non-symmetry of the random walk
does not appear in the appropriate space-time scaling limit.

Throughout the present paper, $O(\cdot)$ stands for the Landau symbol.
When the dependence of 
the $O(\cdot)$ term is
significant, we specify as $O_{N}(\cdot)$, etc.
\section{Framework and Results}
\begin{figure}[tbph]
\begin{center}
\includegraphics[width=8cm]{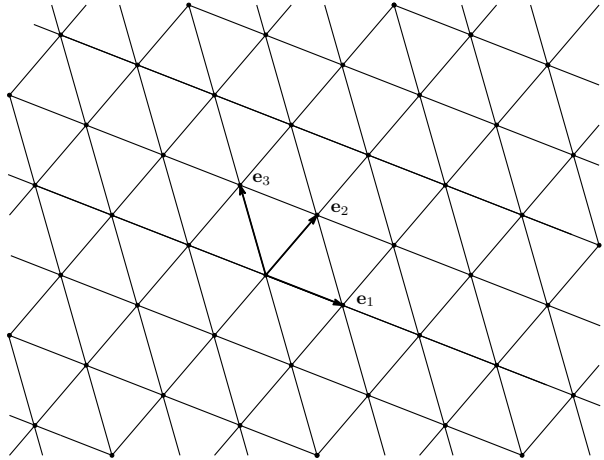} 
\end{center}
\caption{Triangular lattice}
\label{triangular lattice}
\end{figure}

First of all, we prepare some notations and formulate our problem.
Let ${\bf e}_{1}$ and ${\bf e}_{2}$ be linearly independent 
vectors in $\mathbb R^{2}$ and we 
set ${\bf e}_{3}:={\bf e}_{2}-{\bf e}_{1}$ and $K:=\Vert{\bf e}_{1} \Vert_{\mathbb R^{2}}+\Vert{\bf e}_{2} \Vert_{\mathbb R^{2}}$. 
We define the triangular lattice $G=(V,E)$ by
\begin{align*} 
V &=\{ x=x_{1}{\bf e}_{1}+x_{2}{\bf e}_{2} \vert~(x_{1},x_{2}) \in {\mathbb Z}^{2} \}, \\
E &=\left\{ (x,y) \in V\times V \vert~ 
x-y \in \left\{ \pm {\bf e}_1, \pm {\bf e}_2, \pm {\bf e}_3 \right\} \right\}
\end{align*}
(see Figure \ref{triangular lattice}).
We remark that $G$ is isomorphic to
the Cayley graph of the abelian group $\mathbb{Z}^2$ with a set of generators 
$S=\{ \pm (1, 0), \pm (0, 1), \pm (-1, 1) \} $ identified by
$V \simeq \mathbb{Z}^2$ and  $\mathbf{e}_{1} \simeq (1,0)$, 
$\mathbf{e}_{2} \simeq(0,1)$, $\mathbf{e}_{3} \simeq (-1, 1)$.

Throughout the present paper, we impose the following conditions 
on the $1$-step transition probability $p$:
\begin{description}
\item[{\bf (P1):}] 
For every $x\in V$,
\begin{align*}p \big( (x, x+{\bf e}_{1}) \big)&=\alpha, \quad p \big ((x, x-{\bf e}_{1}) \big)=\alpha', \\
p \big((x, x+{\bf e}_{2})\big) &=\beta', \quad  p \big(x,x-{\bf e}_{2}) \big)=\beta, \\
p \big((x, x+{\bf e}_{3}) \big)&=\gamma, \quad
p \big((x, x+{\bf e}_{3}) \big)=\gamma',
\end{align*}
where $\alpha, \alpha', \beta, \beta', \gamma, \gamma' \geq 0$
and $\alpha+\alpha'+\beta+\beta'+\gamma+\gamma'=1$. 
\item[{\bf (P2):}] Zero mean condition:
\begin{equation*}
\sum_{e\in E_{0}}p(e)e=0.
\end{equation*}
\item[{\bf (P3):}] 
\begin{equation*}
\Gamma (p):={\widehat \alpha}{\widehat \beta}+{\widehat \beta}{\widehat \gamma}
+{\widehat \gamma}{\widehat \alpha}>0,
\end{equation*}
where ${\widehat \alpha}:=\alpha+\alpha',
{\widehat \beta}:=\beta+\beta', {\widehat \gamma}:=\gamma+\gamma'$.
\end{description}
%
\begin{re}
Condition {\bf (P2)} is equivalent to the following condition:
\vspace{2mm} \\
$\bullet$
There exists a constant $0\leq \kappa \leq \frac{1}{3}$ 
such that $\alpha-\alpha'=\beta-\beta'=\gamma-\gamma'=\kappa$.
\vspace{2mm} \\
In the case $\kappa=0$, condition {\bf (P2)} implies that $p(e)=p({\overline e})$ holds for every $e\in E$. 
Namely, our random walk is symmetric
with $m_{V}\equiv 1$. On the other hand, by a simple calculation, we see that 
there does not exist the function $m_{V}$ satisfying 
{\rm{(\ref{symmetric})}} unless $\kappa=0$. Hence our random walk is non-symmetric in the case $0< \kappa \leq \frac{1}{3}$.
We can regard the constant $\kappa$ as intensity of the 
non-symmetry. 
\end{re}

We set 
\begin{equation} M_{q}(\theta):=\sum_{e\in E_{0}} p(e) \langle e, \theta \rangle^{q}, 
\quad \theta=
(\theta_{1}, \theta_{2}) \in \mathbb R^{2}, q\in \mathbb N,
\label{Mm}
\end{equation}
where $\langle \cdot, \cdot \rangle$ stands for the scalar product on $\mathbb R^{2}$. Note that 
condition {\bf (P2)} implies $M_{1}(\theta) \equiv 0$. By condition {\bf (P1)}, we also have
\begin{equation}
\vert M_{q}(\theta) \vert \leq \big( \sum_{e\in E_{0}} p(e) \Vert e \Vert_{\mathbb R^{2}}^{q} \big) 
\Vert \theta \Vert_{\mathbb R^{2}}^{q} \leq K^{q} \Vert \theta \Vert_{\mathbb R^{2}}^{q}, 
\quad \theta \in \mathbb R^{2}, q\in \mathbb N.
\label{Mm-est}
\end{equation}

We define the covariance matrix $Q$ 
by 
$$ \langle Q\theta, \theta \rangle =
M_{2}(\theta),
\quad
\theta \in \mathbb R^{2}.$$
In the case of ${\bf e}_{1}={\widehat {\bf e}}_{1}:=\hspace{-1mm}\mbox{ }^{t}(1,0)$
and
${\bf e}_{2}={\widehat {\bf e}}_{2}:=\hspace{-1mm}\mbox{ }^{t}(0,1)$, 
the corresponding covariance matrix is easily calculated as
$$ {\widehat Q}:= \left(
       \begin{array}{cc}
       {\widehat \alpha}+{\widehat \gamma} & -{\widehat \gamma} \\
       -{\widehat \gamma} & {\widehat \beta}+{\widehat \gamma}
       \end{array}
     \right),$$
and hence we obtain ${\rm det}{\widehat Q}=\Gamma(p)$.
For a general pair of two vectors ${\bf e}_{1}$, ${\bf e}_{2}$, we can decompose
the covariance matrix $Q$ as
\begin{equation} Q=T{\widehat Q} \hspace{-1mm}\mbox{ }^{t}T,
\label{Q-decompose}
\end{equation}
where $T=[ {\bf e}_{1}, {\bf e}_{2} ]$ stands for the matrix formed by column vectors ${\bf e}_{1}, {\bf e}_{2}$.
It follows from condition 
{\bf (P3)} and linear independence of ${\bf e}_{1}$, ${\bf e}_{2}$ that 
the covariance matrix $Q$ is positive definite, i.e., 
$\langle \theta, Q\theta \rangle \geq \lambda \Vert \theta \Vert^{2}_{\mathbb R^{2}},
\theta \in \mathbb R^{2}$ for some $\lambda >0$.
By (\ref{Q-decompose}), we also observe that
\begin{equation} \langle Q^{-1}x, y \rangle = \langle {\widehat Q}^{-1}
(x_{1} {\widehat{\bf e}_{1}}+
x_{2} {\widehat{\bf e}_{2}}), y_{1} {\widehat{\bf e}_{1}}+
y_{2} {\widehat{\bf e}_{2}} \rangle 
\label{base-indep}
\end{equation}
holds for all $x=x_{1} {\bf e}_{1}+x_{2} {\bf e}_{2}, y=y_{1}{\bf e}_{1}+y_{2} {\bf e}_{2}\in V$.
This means that the left-hand side of (\ref{base-indep}) is independent of the realization of
the triangular lattice $G$.

Now, we are in a position to state our main result. 
%
Let
\begin{equation}
A(G):= \frac{1}{3 \Gamma(p)^{1/2}}, \quad l := 
\left( \frac{  \widehat{ \beta} + \widehat{ \gamma }  }
{3\Gamma(p)} \right)^{1/2}
\label{constants}
\end{equation}
and we introduce two vectors by
{

$${\mathbf h}_{1}:=\hspace{-1mm}\mbox{ }^{t}(l,0), 
\quad {\mathbf h}_{2}:=\hspace{-1mm}\mbox{ }^{t}
\Big (\frac{\widehat \gamma}{{\widehat \beta}+{\widehat \gamma}} l ,
 \frac{\Gamma(p)^{1/2}}{ \widehat{\beta}+ \widehat \gamma} l \Big)
$$
(see Figure {\ref{standard}}).
\begin{figure}[tbph]
\begin{center}
\includegraphics[width=10cm]{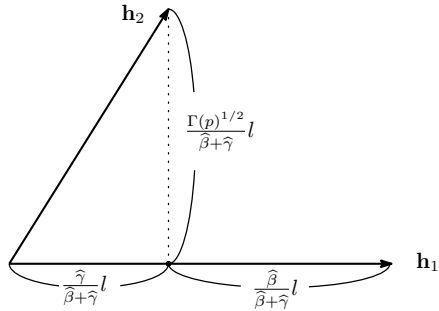} 
\end{center}
\caption{$\mathbf{h}_1$ and $\mathbf{h}_2$}
\label{standard}
\end{figure}
Note that $A(G)$ is equal to 
the area of the parallelogram spanned by $\mathbf{h}_1$ and $\mathbf{h}_2$.
In the case of ${\bf e}_{1}={\mathbf h}_{1}$ and ${\bf e}_{2}={\mathbf h}_{2}$, the corresponding covariance matrix
is 
\begin{equation}
Q= \left( 
\begin{array}{cc}
1/3 & 0 \\
0 & 1/3
\end{array}
\right).
\label{standard-Q}
\end{equation}

Realizing the triangular lattice $G$ along these two vectors,  we have 
\begin{tm}
\label{Main}
{\rm (1)}~
Let us consider the case $0 \leq \kappa < 1/3$. Then we have
\begin{eqnarray*}
2\pi n \cdot 
p(n,x,y) 
&=&
3A(G)
\exp \Big( -\frac{3}{2 n} 
\| y-x \|^2_{\mathbb{R}^2} \Big)
\nonumber 
\\
&\mbox{ }&
\times
\Big( 1+
\sum_{j=1}^{N} n^{-j/2}P_{j}\Big(\frac{y-x}{\sqrt n} 
\Big)
\Big) +O_{N} \big( n^{-\frac{N+1}{2}} \big),
\quad
N\in 
\mathbb N \cup \{0 \},
\end{eqnarray*}
as $n \rightarrow \infty$ uniformly for all $x, y \in V$, where $P_{j}=P_{j}(y)$, 
$j\in \mathbb N$
is a polynomial
of degree 
at most 
$3j$ in the variables $y_{1}, y_{2}$
and an odd or even function depending on whether $j$ is
odd or even.
Furthermore, let us denote the leading term  
of
$\sum_{j=1}^{N} n^{-j/2}P_{j}\big(\frac{y-x}{\sqrt n} \big ), N\geq 2$ 
by $a_1(y-x;\kappa)n^{-1}$, that is, 
\begin{equation*}
\lim_{n \to \infty}
n \Big( 
\sum_{j=1}^{N} n^{-j/2}P_{j}\Big(\frac{y-x}{\sqrt n} \Big )- a_1(y-x ; \kappa) n^{-1}
\Big)=0, \quad x,y\in V.
\end{equation*}
Then the coefficient $a_1(y ;\kappa )$ is explicitly
obtained by
\begin{equation*}
a_1(y;k)=a_1^{(0)}+\kappa a_1^{(1)}(y) +{\kappa}^2 a_1^{(2)},
\end{equation*}
where
\begin{align*}
a_1^{(0)} &= -1+\frac{1}{8\Gamma (p)^2} 
\left\{ \widehat{\alpha} ( \widehat{\beta} +\widehat{\gamma})^2 
+\widehat{\beta} ( \widehat{\gamma} + \widehat{\alpha} )^2
+\widehat{\gamma} ( \widehat{ \alpha} + \widehat{\beta})^2 \right\}, \\
a_1^{(1)}(y)
&= \frac{1}{\Gamma (p)^2}
\left\{ (\widehat{\alpha} \widehat{\beta} -2\widehat{\beta} \widehat{\gamma}
+\widehat{\gamma} \widehat{\alpha} ) y_{1}
+( -\widehat{\alpha} \widehat{\beta} -\widehat{\beta} \widehat{\gamma}
+2\widehat{\gamma}\widehat{\alpha} ) y_{2} \right\}, 
\quad y= y_1 \mathbf{h}_1+y_2 \mathbf{h}_2, \\
a_1^{(2)} &=
\frac{3}{8\Gamma (p)^2 } 
\Big( -1 + \frac{5 \widehat{\alpha} \widehat{\beta} \widehat{\gamma}}
{\Gamma (p)} \Big).
\end{align*}
{\rm (2)}~
Let us consider the case $\kappa=1/3$, i.e., $\alpha=\beta=\gamma=1/3$,  
$\alpha'=\beta'=\gamma'=0$. We set
\begin{equation}
V_{k}:= \big \{ x=x_{1} \widetilde{\mathbf{h}}_1 + x_2 \widetilde{\mathbf{h}}_2 \ | 
\ (x_1, x_2) \in \mathbb{Z}^2 \big \} +k{\mathbf h}_{1}, \quad k=0,1,2,
\label{submodule}
\nonumber
\end{equation}
where $\widetilde{\mathbf{h}}_1=2\mathbf{h}_1 - \mathbf{h}_2$, 
$\widetilde{\mathbf{h}}_2=\mathbf{h}_1 + \mathbf{h}_2$, 
$\mathbf{h}_1= \!^t ( \sqrt{2/3}, 0)$, $\mathbf{h}_2=\!^t ( 1/\sqrt{6}, 1/\sqrt{2} )$
(see Figure {\rm{\ref{periodic}}}). 
For $k,l=0,1,2$, we take $x\in V_{k}$ and $y\in V_{l}$.
Then we have
\begin{eqnarray*}
2\pi n \cdot 
p(n,x,y) 
&=&
9A(G)
\exp \Big( -\frac{3}{2 n} 
\| y-x \|^2_{\mathbb{R}^2} \Big)
\nonumber 
\\
& \mbox{  }&
\times
\Big( 1+
\sum_{j=1}^{N} n^{-j/2}P_{j}\Big(\frac{y-x}{\sqrt{ n}} \Big)
\Big) +O_{N} \big( n^{-\frac{N+1}{2}} \big),
\quad
N\in \mathbb N \cup \{0 \}
\end{eqnarray*}
as $n=3m+(l-k) \rightarrow \infty$ uniformly for all $x$ and $y$. 
In this case, $A(G)=1/{\sqrt 3}$ and the coefficient of the leading term is 
\begin{equation*}
a_1(y-x;1/3)=
-\frac{2}{3}.
\end{equation*}
\end{tm}
\begin{figure}[tbph]
\begin{center}
\includegraphics[width=8cm]{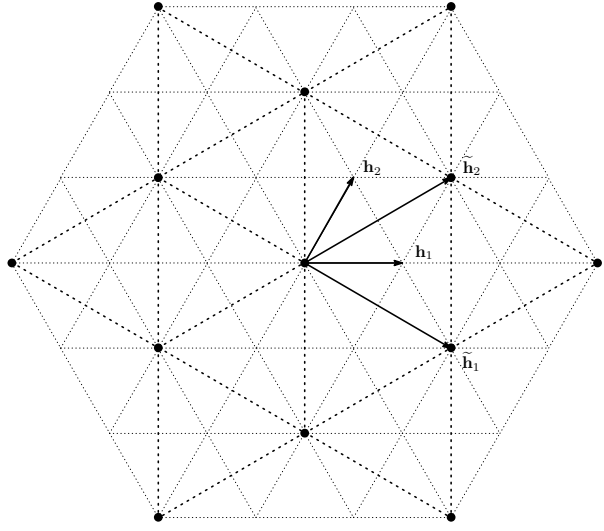}
\end{center}
\caption{$V_0$ in the case $\kappa=1/3$}
\label{periodic}
\end{figure}

We characterize the pair $\mathbf{h}_1$, $\mathbf{h}_2$ 
through a variational problem on the crystal lattices based on the idea 
of Kotani and Sunada
\cite{KS00, KS00-2}.

 We set ${\bf e}_{1}:= \!^t(u,0), {\bf e}_{2}:=\!^t (v_1,v_2)$.
Without loss of generality, we may assume 
$u, v_2>0$.
Then 
the energy of the (quotient graph of the) triangular lattice  $G$  is given by
\begin{eqnarray}
{\cal E}(G)
&:=& \frac{1}{2} \sum_{e \in E_{0}} p(e) \| t(e)-o(e) \|^{2}_{\mathbb R^{2}}
\nonumber \\
&=& \frac{1}{2} \Big [ {\widehat \alpha}u^2 +{\widehat \beta}(v_1^{2}+v_2^{2})+
{\widehat \gamma} \big \{ 
(v_1-u)^{2}+v_2^{2} \big \} \Big ].
\nonumber
\end{eqnarray}
Minimizing ${\cal E}(G)$ with respect to $(u, v_1, v_2)$ under the condition
${\rm det}
( \langle {\bf e}_{i}, {\bf e}_{j}
\rangle )_{i,j=1}^{2}=A(G)$,
we obtain 
$$ u=l,  \quad
v_{1}=\frac{ \widehat \beta}{\widehat \beta +\widehat \gamma} l, 
\quad v_{2}= \frac{ \Gamma(p)^{1/2}}{ \widehat{ \beta} + \widehat{\gamma} } l.
$$
Hence we have derived $\mathbf{e}_1=\mathbf{h}_1$, $
\mathbf{e}_2= \mathbf{h}_2$.
\begin{re}
\label{Traiangular explain}
In the case 
$
\alpha=\alpha^\prime=\beta=\beta^\prime=\gamma=\gamma^\prime=1/6,
$
the random walk is said to be simple. 
 In this case, we have
\begin{equation*}
A(G)= \frac{1}{\sqrt{3}}, \quad l=\sqrt{ \frac{2}{3} }, \quad 
\mathbf{h}_1= \!^t \left( \sqrt{ \frac{2}{3} }, 0\right), \quad 
\mathbf{h}_2= \!^t \left( \frac{1}{\sqrt{6}}, \frac{1}{\sqrt{2}} \right).
\end{equation*}
We note the volume of the Albanese torus is equal to $1/\sqrt{3}$
 (see {\rm{\cite[page 640]{KS00}}}, {\rm{\cite[Section 10]{KSS}}} ).  
The standard realization of the triangular lattice 
is the equilateral triangular lattice in $\mathbb{R}^2$ each of whose edge has length 
$\sqrt{2/3}$.
The quantity $A(G)$ introduced in {\rm{(\ref{constants})}} can be regarded 
as a generalization of the notion of the 
volume of the Albanese torus to some non-symmetric cases. 
Henceforth, we call the realization of the triangular lattice $G$ along 
$\mathbf{h}_1$ and  $\mathbf{h}_2$ the standard realization even if the 
random walk is non-symmetric.
\end{re}

Next, let us consider the {\it discrete heat equation} on the triangular lattice $G=(V,E)$
with initial condition $u_{0}:V\to \mathbb R$:
\begin{equation}
\left\{
\begin{array}{rl}
  \left( \partial_1 +\Delta_d \right) u(n,x)&=0, \qquad (n,x)\in {\mathbb Z}_{+} \times V,\\
  u(0,x) &=u_0(x), \qquad x\in V.
 \end{array}
\right. 
\label{heat}
\end{equation}
Here $\partial_1 u(n,x):=u(n+1,x)-u(n,x)$ and the {\it discrete Laplacian} $\Delta_d$
is defined by
\begin{equation*}
\Delta_d u(n,x):= \sum_{e \in E_x} p(e) \Big( u(n,o(e))- u(n,t(e)) \Big).
\end{equation*}
The operator $L:=I-\Delta_d$ is called the {\it transition operator} 
associated with the random walk on $G$. It is easy to see that 
the solution of (\ref{heat}) is given by
\begin{equation*}
u(n,x)=L^n u_0(x).
\end{equation*}
We note that the solution is rewritten by
\begin{equation*}
u(n,x)=\sum_{y \in V} p(n,x,y)u_0(y).
\end{equation*}

For $t>0$, let $H_t$ be the heat operator defined by
\begin{equation*}
H_t f(x)=\frac{1}{2\pi t   }\int_{\mathbb{R}^2} 
\exp \Big(
- \frac{\| x- z\|_{\mathbb{R}^2}^2 }{2t} \Big )
f(z)dz, \quad f\in C_{\infty}(\mathbb R^{2}),
\end{equation*}
where $C_{\infty}(\mathbb R^{2})$ is the set of continuous functions on $\mathbb R^{2}$ vanishing
at infinity.
Note that $\{ H_t \}_{t\geq 0}$ is a semigroup 
whose infinitesimal generator is 
\begin{equation*}
\frac{1}{2}\Delta:=\frac{1}{2}
\left( \frac{\partial^2}{\partial x_1^2} + \frac{\partial^2}{\partial x_2^2} \right).
\end{equation*}

Now, we take the standard realization of the triangular lattice $G$ as before.
Applying Theorem \ref{Main}, we obtain the following.
\begin{tm} \label{FTCLT}
Let $t>0$ and 
$\{\delta_n \}_{n=1}^{\infty}$
be a sequence of positive real numbers satisfying $\lim_{n\to \infty} n\delta_{n}^2=  3t$.
Then for every continuous function $f$ on $\mathbb R^{2}$ with compact support
and for a sequence $\{ x_n \}_{n=1}^{\infty}$ in $V$ with 
$\lim_{n \rightarrow \infty} \delta_n x_n =x \in \mathbb{R}^2$,
we have
\begin{equation*}
\lim_{n \rightarrow \infty}
L^n (f\circ \delta_n) (x_n)= H_t f(x).
\end{equation*}
\end{tm}

This theorem is also an immediate consequence of the following result via approximation 
theory due to Trotter \cite{Trotter} (see also Kotani \cite{Kotani}).
\begin{tm}\label{Trotter}
\begin{enumerate}
\item[{\rm{(1)}}]
\begin{equation*}
\lim_{\delta \rightarrow 0} \left\| \frac{3}{ \delta ^2} \Delta_d (f\circ \delta) 
+\frac{1}{2}(\Delta f) \circ \delta \right\|_\infty = 0, \quad f\in C^{\infty}_{0}(\mathbb R^{2})
\end{equation*}
\item[{\rm{(2)}}]Let $\{ \delta_n \}_{n=1}^\infty$ be a sequence as in Theorem {\rm{\ref{FTCLT}}}. 
\begin{equation*}
\lim_{n\rightarrow \infty} 
\| L^n (f\circ \delta_n) -(H_t  f) \circ \delta_n \|_\infty = 0, \quad 
 f \in C_\infty (\mathbb{R}^2 ).
\end{equation*}
\end{enumerate}
\end{tm}
%
\section{Preliminaries}
In this section, we give some basic facts for the proof of
Theorem \ref{Main}. In what follows, we denote
${\bf h}_{3}:={\bf h}_{2}-{\bf h}_{1}$,
$x_{3}:=x_{2}-x_{1}$
and $(\frac{\partial}{\partial x_{3}}):=
(\frac{\partial}{\partial x_{2}})-(\frac{\partial}{\partial x_{1}})$ for convenience.

We begin with the explicit form of $M_{q}(\theta), q\in \mathbb N$ defined in (\ref{Mm}).
\begin{lm}
We have
\begin{equation}
M_{q}(\theta)=
\begin{cases} {\kappa} \big( \langle {\bf e}_{1}, \theta \rangle^{q}-
\langle {\bf e}_{2}, \theta \rangle^{q}
+\langle  {\bf e}_{3}, \theta \rangle^{q} \big )
& \text{ (if $q$ is odd) },
\\
{\widehat \alpha}
\langle {\bf e}_{1}, \theta \rangle^{q}
+
{\widehat \beta}
\langle {\bf e}_{2}, \theta \rangle^{q}
+
{\widehat \gamma}
\langle  {\bf e}_{3}, \theta \rangle^{q} 
& \text{ (if $q$ is even)}.
\end{cases}
\nonumber
\end{equation}
In particular, $M_{q}(\theta)\equiv 0$ for every odd number $q$ in the
case where the random walk is symmetric.
\end{lm}

Noting (\ref{base-indep}) and (\ref{standard-Q}), we easily see 
$
\langle Q^{-1}{\bf e}_{i}, {\bf e}_{j} \rangle =3\langle {\bf h}_{i}, {\bf h}_{j} \rangle,
i=1,2,3
$. A direct calculation of the inner product $\langle {\bf h}_{i}, {\bf h}_{j} \rangle, 
i,j=1,2,3$
yields the following lemma:
\begin{lm}
\label{naiseki-keisan}
For any pair of linearly independent vectors ${\bf e}_{1}, {\bf e}_{2}\in \mathbb R^{2}$, we have
\begin{eqnarray}
\langle Q^{-1}{\bf e}_{1}, {\bf e}_{1} \rangle &=&
\frac{ \widehat \beta +\widehat \gamma}{\Gamma(p)},
\quad
\langle Q^{-1}{\bf e}_{2}, {\bf e}_{2} \rangle 
=
\frac{ \widehat \gamma+ \widehat \alpha}{\Gamma(p)},
\quad
\langle Q^{-1}{\bf e}_{3}, {\bf e}_{3} \rangle 
=\frac{ \widehat \alpha +\widehat \beta}{\Gamma(p)},
\nonumber
\\
\langle Q^{-1}{\bf e}_{1}, {\bf e}_{2} \rangle &=&
\frac{ \widehat \gamma }{\Gamma (p)},
\qquad
\langle Q^{-1}{\bf e}_{2}, {\bf e}_{3} \rangle =
\frac{ \widehat \alpha}{\Gamma(p)},
\qquad
\langle Q^{-1}{\bf e}_{3}, {\bf e}_{1} \rangle =
\frac{ -\widehat \beta}{\Gamma(p)}.
\nonumber
\end{eqnarray}
\end{lm}

Next, we recall an elementary fact about the Fourier transform
\begin{equation}
\int_{\mathbb R^{2}} 
e^{-\frac{1}{2} \langle Q \theta, \theta \rangle} 
e^{-{\sqrt{-1}} \langle x , \theta \rangle} d\theta
=2\pi ({\rm det}Q)^{-1/2} e^{-\frac{1}{2} \langle Q^{-1} x, x \rangle}, \quad x\in \mathbb R^{2}.
\label{Gauss-Fourier}
\end{equation}
Differentiating both sides of (\ref{Gauss-Fourier}) with respect to $x_{i},~i=1,2,3$, we have
\begin{eqnarray}
\lefteqn{
\int_{\mathbb R^{2}} 
\langle {\bf e}_{i}, \theta \rangle
e^{-\frac{1}{2} \langle Q \theta, \theta \rangle} 
e^{-{\sqrt{-1}} \langle x , \theta \rangle} d\theta}
\nonumber \\
&=&
2\pi {\sqrt{ -1}} 
\cdot({\rm det}Q)^{-1/2} 
\big( -\langle Q^{-1} {\bf e}_{i}, x \rangle \big )
e^{-\frac{1}{2} \langle Q^{-1} x, x \rangle}, \quad x\in \mathbb R^{2}, i=1,2,3.
\nonumber
\end{eqnarray}
Repeating this argument 
several times,
we obtain the following proposition:
\begin{pr}
\label{3.3}
Let us set 
$$ F(i_{1}, \ldots, i_{N})(x):=\int_{\mathbb R^{2}} \big(
\prod_{k=1}^{N} \langle {\bf e}_{i_{k}}, \theta \rangle \big) 
e^{-\frac{1}{2} \langle Q \theta, \theta \rangle} 
e^{-{\sqrt{-1}} \langle x , \theta \rangle}
d\theta, \quad N\in \mathbb N,~ i_{1}, \ldots, i_{N}=1,2,3.
$$ 
Then we have
\begin{equation} F(i_{1}, \ldots, i_{N})(x)=2\pi({\sqrt{-1}})^{N} ({\rm det}Q)^{-1/2}
e^{-\frac{1}{2} \langle Q^{-1} x, x \rangle}
G(i_{1},\ldots, i_{N})(x),
\nonumber
\end{equation}
where $\{G(i_{1},\ldots, i_{k})(x): x\in \mathbb R^{2}, k=1,\ldots, N \}$ is determined as the
solution of the recursive system of the following equations starting from $k=1$ to $k=N$:
\begin{equation}
\left\{ \begin{array}{ll}
G(i_{1},\ldots, i_{k})(x)
= - \langle Q^{-1}{\bf e}_{i_{k}}, x \rangle G(i_{1},\ldots, i_{k-1})(x)
+{\displaystyle{\frac{\partial}{\partial x_{i_{k}}} G(i_{1},\ldots, i_{k-1})(x)}},
\\ \\
G(i_{1})(x)=-\langle Q^{-1}{\bf e}_{i_{1}}, x \rangle.
\end{array} \right. \label{Recursive}
\end{equation}
\end{pr}
%
\begin{re}
\label{3.4}
$G(i_{1},\ldots, i_{N})(x)$ is decomposed by
\begin{equation*}
G(i_{1},\ldots, i_{N})(x)=
\begin{cases}
\displaystyle{\sum_{l=0}^{\frac{N-1}{2}} G(i_{1},\ldots, i_{N})_{(2l+1)}(x)}
& \text{ (if $N$ is odd)},
\\ 
\displaystyle{\sum_{l=0}^{\frac{N}{2}} G(i_{1},\ldots, i_{N})_{(2l)}(x) }
& \text{ (if $N$ is even) },
\end{cases}
\end{equation*}
where
$G(i_{1},\ldots, i_{N})_{(k)}(x)$, $k=0,1,\ldots, N$
is a homogeneous polynomial 
of degree $k$
in the variables $\langle Q^{-1}{\bf e}_{1},x \rangle$,
$\langle Q^{-1}{\bf e}_{2},x \rangle$ and $\langle Q^{-1}{\bf e}_{3},x \rangle$.
In particular, 
\begin{eqnarray}
G(i_{1},i_{2})(x)&=& \langle Q^{-1}{\bf e}_{i_1},x \rangle
\langle Q^{-1}{\bf e}_{i_2},x \rangle -\langle Q^{-1}{\bf e}_{i_1},{\bf e}_{i_2}\rangle,
\qquad i_{1},i_{2}=1,2,3,
\nonumber
\\
G(i_{1},i_{2},i_{3})(x)&=&-\langle Q^{-1}{\bf e}_{i_1},x \rangle
\langle Q^{-1}{\bf e}_{i_2},x \rangle
\langle Q^{-1}{\bf e}_{i_3},x \rangle 
\nonumber \\
&\mbox{ }&
+\langle Q^{-1}{\bf e}_{i_1},{\bf e}_{i_2} \rangle
\langle Q^{-1}{\bf e}_{i_3},x \rangle
+\langle Q^{-1}{\bf e}_{i_2},{\bf e}_{i_3} \rangle
\langle Q^{-1}{\bf e}_{i_1},x \rangle
\nonumber \\
&\mbox{ }&
+\langle Q^{-1}{\bf e}_{i_3},{\bf e}_{i_1} \rangle
\langle Q^{-1}{\bf e}_{i_2},x \rangle,
\qquad i_{1},i_{2},i_{3}=1,2,3.
\nonumber
\end{eqnarray}
\end{re}

We define the characteristic function of the $1$-step transition probability $p$ by
\begin{equation}
\varphi (\theta):= \sum_{e\in E_{0}} p(e) \exp \{ {\sqrt{-1}} \langle e, \theta \rangle \}, \quad \theta=(\theta_{1},\theta_{2})
\in \mathbb R^{2}.
\nonumber
\end{equation}
We denote $\varphi$ especially by $\phi$ in the case of
${\bf e}_{1}={\widehat {\bf e}}_{1}$ and 
${\bf e}_{2}={\widehat {\bf e}}_{2}$. 
Noting that the characteristic function of the $n$-step transition probability of the random walk 
starting at the origin is equal to $\varphi^{n}$,
we have the following integral expression of $p(n,x,y)$.
\begin{lm}\label{Fourier}
\begin{equation*}
p(n,x,y)=\frac{ \vert {\rm det}T \vert}{(2\pi)^{2}}
\int_{(\hspace{-1mm}\mbox{ }^{t}T)^{-1}(D)}
\varphi (\theta)^{n} \exp \{ -\sqrt{-1} \langle y-x, \theta \rangle \}
d\theta, \quad x,y \in V,~ n\in \mathbb N,
\end{equation*}
where $D:=[-\pi, \pi]^{2}$.
\end{lm}
{\bf{Proof.}}~Since $p(n,x,y)=p(n,0,y-x)$, we may suppose $x=0$.
We also note that 
$p(n,0,y)$, $y=y_{1}{{\bf e}}_{1}+y_{2}{{\bf e}}_{2}$
is independent 
of the pair of linearly independent two vectors ${\bf e}_{1}, {\bf e}_{2}$.
In the case ${\bf e}_{1}={\widehat {\bf e}}_{1}$, ${\bf e}_{2}={\widehat {\bf e}}_{2}$, 
the following identity is well-known:
\begin{equation}
p(n,0,y)=\frac{1}{(2\pi)^{2}}
\int_{D}
\phi (\theta)^{n} \exp \big \{ -\sqrt{-1} 
\sum_{i=1}^{2}y_{i}\theta_{i}
\big \}d\theta,
\quad
y=y_{1}{\widehat {\bf e}}_{1}+y_{2}{\widehat {\bf e}}_{2},~ n\in \mathbb N.
\label{well-known}
\end{equation}
See Lawler and Limic \cite[Section 2.2.2]{LL} for details.

For a general pair of independent two vectors ${{\bf e}}_{1}, {{\bf e}}_{2}$, 
we observe
\begin{equation}
\varphi (\theta)=\phi (\hspace{-1mm}\mbox{ }^{t}T \theta), \quad \theta \in \mathbb R^{2}.
\label{phi-varphi}
\end{equation}
Then it follows from (\ref{well-known}) and (\ref{phi-varphi}) that
\begin{eqnarray}
\lefteqn{
\frac{1}{(2\pi)^{2}}
\int_{D}
\phi (\theta)^{n} \exp \big \{ -\sqrt{-1} 
\sum_{i=1}^{2}y_{i}\theta_{i}
\big \}d\theta}
\nonumber \\
&=&
\frac{1}{(2\pi)^{2}}
\int_{D}
\varphi \big( (\hspace{-1mm}\mbox{ }^{t}T)^{-1}\theta \big)^{n} \exp \big \{ -\sqrt{-1} 
\big \langle \sum_{i=1}^{2}y_{i}T{\widehat {\bf e}}_{i}, (\hspace{-1mm}\mbox{ }^{t}T)^{-1}\theta \big \rangle
\big \}d\theta
\nonumber \\
&=&\frac{1}{(2\pi)^{2}}
\int_{D}
\varphi \big( (\hspace{-1mm}\mbox{ }^{t}T)^{-1}\theta \big)^{n} \exp \big \{ -\sqrt{-1} 
\langle y, (\hspace{-1mm}\mbox{ }^{t}T)^{-1}\theta \rangle
\big \}d\theta
\nonumber \\
&=&
\frac{ \vert {\rm det}T \vert}{(2\pi)^{2}}
\int_{(\hspace{-1mm}\mbox{ }^{t}T)^{-1}(D)}
\varphi (\theta')^{n} \exp \{ -\sqrt{-1} \langle y, \theta' \rangle \}
d\theta', \quad 
y=y_{1}{ {\bf e}}_{1}+y_{2}{{\bf e}}_{2},~n\in \mathbb N,
\nonumber
\end{eqnarray}
where we performed the change of variables
$\theta=\hspace{-1mm}\mbox{ }^{t}T \theta'$ for the final line.
This completes the proof.
\qed
\begin{lm}
\label{AP-test}
{\rm (1)}~In the case $0\leq \kappa<1/3$, 
the characteristic function $\varphi(\theta)$ defined on $(\hspace{-1mm}\mbox{ }^{t}T)^{-1}(D)$
has the following property:
$\vert \varphi (\theta) \vert=1$
holds only when $\theta=0$.
\\
{\rm (2)}~In the case $\kappa=1/3$, the characteristic function $\varphi(\theta)$
does not satisfy the above property. 
\end{lm}
{\bf Proof.}~First, we prove (1).
By (\ref{phi-varphi}),
it is sufficient to show that $\vert \phi (\theta) \vert=1, \theta \in D$ implies $\theta=0$.
We calculate the characteristic function $\phi(\theta)$ as
\begin{eqnarray}
\phi(\theta)&=& \big( {\widehat \alpha} \cos (\theta_{1})+\sqrt{-1} k\sin(\theta_{1}) \big)
+\big( {\widehat \beta} \cos (-\theta_{2})+\sqrt{-1} k\sin(-\theta_{2}) \big)
\nonumber \\
&\mbox{ }&+\big( {\widehat \gamma} \cos (\theta_{2}-\theta_{1})+\sqrt{-1} k\sin(\theta_{2}-\theta_{1}) \big)
\nonumber \\
&=:&\phi_{1}(\theta)+\phi_{2}(\theta)+\phi_{3}(\theta),
\quad \theta=(\theta_{1}, \theta_{2}) \in D.
\label{sin-cos}
\end{eqnarray}
Note that the assumption $0\leq \kappa <\frac{1}{3}$ implies 
$\min \{ {\widehat \alpha}, {\widehat \beta}, {\widehat \gamma} \}> \kappa \geq 0$. Then we have
\begin{equation}
\vert \phi_{1}(\theta) \vert \leq {\widehat \alpha},
\quad
\vert \phi_{2}(\theta) \vert \leq {\widehat \beta},
\quad
\vert \phi_{3}(\theta) \vert \leq {\widehat \gamma}, \quad \theta=(\theta_{1}, \theta_{2}) \in D.
\label{daen}
\end{equation} 
We also observe that
and 
$\vert \phi_{1}(\theta)\vert={\widehat \alpha},
\vert \phi_{2}(\theta)\vert= {\widehat \beta}$ and
$\vert \phi_{3}(\theta) \vert=
{\widehat \gamma}$ imply
$\phi_{1}(\theta)=\pm {\widehat \alpha}$, 
$\phi_{2}(\theta)=\pm {\widehat \beta}$ and $\phi_{3}(\theta)=\pm {\widehat \gamma}$, respectively.

Now, we suppose $\vert \phi (\theta) \vert=1$ on $D$. By combining (\ref{sin-cos}) and (\ref{daen})
with ${\widehat \alpha}+{\widehat \beta}+{\widehat \gamma}=1$, we deduce 
$$(\phi_{1}(\theta), \phi_{2}(\theta), \phi_{3}
(\theta))=({\widehat \alpha}, {\widehat \beta}, {\widehat \gamma}), (-{\widehat \alpha}, -{\widehat \beta}, -{\widehat \gamma}),
\quad \theta \in D.$$
We easily see that the first equation 
$(\phi_{1}(\theta), \phi_{2}(\theta), \phi_{3}
(\theta))=({\widehat \alpha}, {\widehat \beta}, {\widehat \gamma})$ 
has the solution
$(\theta_{1}, \theta_{2})=(0,0)$. On the other hand, there exists
no solution 
of the second equation
$(\phi_{1}(\theta), \phi_{2}(\theta), \phi_{3}
(\theta))=(-{\widehat \alpha}, -{\widehat \beta}, -{\widehat \gamma})$. Hence we 
conclude $\theta=0$, which completes the proof of (1).

The item (2) is obvious by recalling $\alpha=\beta=\gamma=1/3$ and 
$\alpha'=\beta'=\gamma'=0$. Actually, $\vert \phi (\theta) \vert=1$ 
on $D$ implies $(\theta_{1},\theta_{2})=(0,0), (\frac{2\pi}{3}, -\frac{2\pi}{3})$, $(-\frac{2\pi}{3}, \frac{2\pi}{3})$.
\qed
\begin{re}In the case $0\leq \kappa<1/3$, the random walk is aperiodic, that is, 
the period of the random walk $d(p):={\rm gcd} \{ n\in \mathbb N:
p(n,0,0) >0 \}$ is equal to $1$. It is derived from 
{\rm (1)} of Lemma {\rm{\ref{AP-test}}}. See e.g., Spitzer {\rm{\cite[P8 in Section 7]{Spitzer}}}
for the proof.
On the other hand, 
the random walk is periodic with
$d(p)=3$ in the case $\kappa=1/3$. 
\end{re}
%

Before closing this section, we present an asymptotic expansion formula of
the characteristic function $\varphi$ 
which plays a crucial role in the next section. 
We set 
\begin{equation}
\chi_{q}(\theta):=({\sqrt{-1}})^{-q} \big ( \frac{d}{dt} \big)^{q} \Big \vert_{t=0}
\log \varphi(t\theta), \quad \theta\in \mathbb R^{2}, q\in \mathbb N.
\label{cum}
\end{equation}
(See Bhattacharya and Ranga Rao \cite[page 47]{BR}.)
We note that $\chi_{q}(\theta)$ is a polynomial in the variables $M_{2}(\theta), \ldots, M_{q}(\theta)$.
In particular, $\chi_{1}(\theta) \equiv 0, \chi_{2}(\theta)=M_{2}(\theta), \chi_{3}(\theta)=M_{3}(\theta)$ and 
$\chi_{4}(\theta)=M_{4}(\theta)-3M_{2}(\theta)^{2}$.

The following proposition is taken from \cite[Lemma 7.1 and Theorem 9.11]{BR}.
\begin{pr} \label{phi-Taylor}
Let $n \in \mathbb N$ and $N\in \mathbb N \cup \{0 \}$. 
Then there exist 
positive 
constants $C_{1}(N), C_{2}(N)$ such that for all 
$\theta \in \mathbb R^{2}$ with
\begin{equation}  \langle Q\theta, \theta \rangle^{1/2}
\leq C_{1}(N) \Big ( \sum_{e\in E_{0}} p(e) \langle Q^{-1}e,e \rangle^{N/2} \Big )^{-\frac{1}{N+3}} \cdot n^{\frac{N+1}{2(N+3)}}, 
\label{check}
\nonumber
\end{equation}
we have
\begin{eqnarray}
\lefteqn{ \Big \vert
\varphi \Big (\frac{\theta}{\sqrt n} \Big)^{n}- \exp \Big(-\frac{1}{2} \langle Q\theta, \theta \rangle \Big)
\cdot
\Big ( b_{0}(\theta)+b_{1}(\theta)n^{-1/2}+\cdots +b_{N}(\theta)n^{-N/2} \Big )
\Big \vert}
\nonumber \\
& &\leq 
C_{2}(N) \exp \Big(-\frac{1}{4} \langle Q\theta, \theta \rangle \Big) 
\Big( \langle Q\theta, \theta \rangle^{N+3}+\langle Q\theta, \theta \rangle^{3(N+1)} \Big) \cdot n^{-\frac{N+1}{2}}.
\nonumber
\end{eqnarray}
Here $b_{0}(\theta)\equiv 1$ and 
$b_{j}(\theta)$, $j=1,\ldots, N$ is 
written as
\begin{equation}
b_{j}(\theta)=({\sqrt -1})^{3j}b_{j}^{(3j)}(\theta)+
({\sqrt -1})^{3j-2}b_{j}^{(3j-2)}(\theta)+\cdots +
({\sqrt -1})^{j+2}b_{j}^{(j+2)}(\theta), 
\label{BJ}
\end{equation}
where $b_{j}^{(k)}(\theta)$, $k=j+2,j+4,\ldots, 3j$ is a 
polynomial in the variables $\chi_{3}(\theta), \ldots, \chi_{k}(\theta)$
and it can be regarded  as a
homogeneous
polynomial of degree $k$
in the variables $\langle {\bf e}_{1}, \theta \rangle, \langle {\bf e}_{2}, \theta \rangle$ and
$\langle {\bf e}_{3}, \theta \rangle$.
\end{pr}
\begin{re}
\label{R-3-9}
In particular, 
\begin{eqnarray}
b_{1}(\theta)&=&
({\sqrt -1})^{3} \Big (\frac{\chi_{3}(\theta)}{6}\Big)
=({\sqrt -1})^{3} \Big (\frac{M_{3}(\theta)}{6}\Big), 
\nonumber
\\
b_{2}(\theta)&=&
({\sqrt -1})^{6} \Big(\frac{\chi_{3}(\theta)^{2}}{72} \Big)+ ({\sqrt -1})^{4}
\Big( \frac{\chi_{4}(\theta)}{24} \Big)
\nonumber \\
&=&({\sqrt -1})^{6} \Big(\frac{M_{3}(\theta)^{2}}{72} \Big)+ ({\sqrt -1})^{4}
\Big( \frac{M_{4}(\theta)}{24}-\frac{M_{2}(\theta)^{2}}{8} \Big).
\nonumber
\end{eqnarray}
\end{re}
\section{Proof of the theorems}
\subsection {Proof of Theorem \ref{Main}}
In this subsection, we prove Theorem \ref{Main2} below.
Note that we easily obtain Theorem \ref{Main} by combining (\ref{standard-Q}) and  
Lemma \ref{naiseki-keisan} with Theorem \ref{Main2}.
\begin{tm}
\label{Main2}
{\rm (1)}~
Let us consider the case $0 \leq \kappa < 1/3$. Then we have
\begin{eqnarray}
2\pi n \cdot 
p(n,x,y) 
&=&
3A(G)
\exp \Big( -\frac{1}{2n} \langle Q^{-1}(y-x), y-x \rangle \Big)
\nonumber 
\\
& \mbox{  }&
\times
\Big( 1+
\sum_{j=1}^{N} n^{-j/2}P_{j}\big(\frac{y-x}{\sqrt n} \big)
\Big) +O_{N} \big( n^{-\frac{N+1}{2}} \big),
\quad
N \in \mathbb N \cup \{0 \}
\label{asymptotic-general}
\end{eqnarray}
as $n \rightarrow \infty$ uniformly for all $x, y \in V$, where $P_{j}=P_{j}(y)$,
$j\in \mathbb N$
is a polynomial of degree 
at most 
$3j$ in the variables $y_{1},y_{2}$ and an odd or even function depending on whether $j$ is
odd or even. Furthermore, let us denote the leading term  
of
$\sum_{j=1}^{N} n^{-j/2}P_{j}\big(\frac{y-x}{\sqrt n} \big ), N\geq 2$ 
by $a_1(y-x;\kappa)n^{-1}$.
Then the coefficient $a_1(y ;\kappa )$ is explicitly
obtained by
\begin{equation}
a_{1}(y;\kappa)= a^{(0)}_{1}+\kappa a^{(1)}_{1}(y)+{\kappa}^{2} a^{(2)}_{1},
\label{explicit}
\end{equation}
where
\begin{eqnarray}
a_{1}^{(0)}
&=& -1+
\frac{1}{8} \Big ( {\widehat \alpha} \langle Q^{-1}{\bf e}_{1}, {\bf e}_{1} \rangle^{2}
+{\widehat \beta} \langle Q^{-1}{\bf e}_{2}, {\bf e}_{2} \rangle^{2}
+
{\widehat \gamma} \langle Q^{-1}{\bf e}_{3},{\bf e}_{3} \rangle^{2}
\Big ),
\nonumber \\
a_{1}^{(1)}(y)&=&\frac{1}{2}
\Big \{ \langle Q^{-1}{\bf e}_{1}, {\bf e}_{2}
\rangle
\langle Q^{-1}{\bf e}_{3}, y \rangle
+\langle Q^{-1}{\bf e}_{2}, {\bf e}_{3} \rangle
\langle Q^{-1}{\bf e}_{1},y \rangle
+
\langle Q^{-1}{\bf e}_{3}, {\bf e}_{1} \rangle
\langle Q^{-1}{\bf e}_{2},y \rangle
\Big \},
\nonumber \\
a_{1}^{(2)}&=& -\frac{5}{24} \Big ( \langle Q^{-1}{\bf e}_{1}, {\bf e}_{1} \rangle^{3}+
\langle Q^{-1} {\bf e}_{2}, {\bf e}_{2} \rangle^{3}
+
\langle Q^{-1}{\bf e}_{3}, {\bf e}_{3} \rangle^{3}
\Big )
\nonumber \\
&\mbox{ }&
+\frac{1}{6}
\Big ( \langle Q^{-1}{\bf e}_{1}, {\bf e}_{2} \rangle^{3}+
\langle Q^{-1}{\bf e}_{2}, {\bf e}_{3} \rangle^{3}
-
\langle Q^{-1}{\bf e}_{3}, {\bf e}_{1} \rangle^{3}
\Big )
\nonumber \\
&\mbox{ }&
+\frac{1}{4}
\Big (
 \langle Q^{-1}{\bf e}_{1}, {\bf e}_{1} \rangle
\langle Q^{-1}{\bf e}_{1}, {\bf e}_{2} \rangle
\langle Q^{-1}{\bf e}_{2}, {\bf e}_{2} \rangle
+
\langle Q^{-1}{\bf e}_{2}, {\bf e}_{2} \rangle
\langle Q^{-1}{\bf e}_{2}, {\bf e}_{3} \rangle
\langle Q^{-1}{\bf e}_{3}, {\bf e}_{3} \rangle
\nonumber \\
&\mbox{ }&~~~~~
-\langle Q^{-1}{\bf e}_{3}, {\bf e}_{3} \rangle
\langle Q^{-1}{\bf e}_{3}, {\bf e}_{1} \rangle
\langle Q^{-1}{\bf e}_{1}, {\bf e}_{1} \rangle
\Big ).
\label{general-a1}
\nonumber
\end{eqnarray}
{\rm (2)}~Let us consider the case $\kappa=1/3$, i.e., $\alpha=\beta=\gamma=1/3$,  
$\alpha'=\beta'=\gamma'=0$. We set
\begin{equation}
V_{k}:= \big \{ x=x_{1} \widetilde{\mathbf{e}}_1 + x_2 \widetilde{\mathbf{e}}_2 \ | 
\ (x_1, x_2) \in \mathbb{Z}^2 \big \} +k{\mathbf e}_{1}, \quad k=0,1,2,
\label{submodule}
\end{equation}
where $\widetilde{\mathbf{e}}_1=2\mathbf{e}_1 - \mathbf{e}_2$, 
$\widetilde{\mathbf{e}}_2=\mathbf{e}_1 + \mathbf{e}_2$.
For $k,l=0,1,2$, we take $x\in V_{k}$ and $y\in V_{l}$.
Then we have
\begin{eqnarray}
2\pi n \cdot 
p(n,x,y) 
&=&
9A(G)
\exp \Big( -\frac{1}{2n} \langle Q^{-1}(y-x), y-x \rangle \Big)
\nonumber 
\\
& \mbox{  }&
\times
\Big( 1+
\sum_{j=1}^{N} n^{-j/2}{P}_{j}\Big(\frac{y-x}{\sqrt{ n}} \Big)
\Big) +O_{N} \big( n^{-\frac{N+1}{2}} \big),
\quad
N\in \mathbb N \cup \{0 \}
\label{asymptotic-3}
\end{eqnarray}
as $n=3m+(l-k) \rightarrow \infty$ uniformly for all $x$ and $y$.
In this case, the coefficient of the
leading term is 
\begin{equation*}
a_1(y-x;1/3)=
-\frac{2}{3}.
\end{equation*}
\end{tm}
{\bf Proof.} As we mentioned in the proof of Lemma \ref{Fourier},
we may suppose $x=0$ throughout the proof. 

First of all, we prove (1).
By recalling that the covariance matrix $Q$ is positive definite and 
Lemma \ref{AP-test}, 
we can choose a positive constant $R$ sufficiently small such that both
\begin{equation}
\label{outside-est}
\eta :=\sup \{ \vert \varphi (\theta) \vert: \theta \in (\hspace{-1mm}\mbox{ }^{t}T)^{-1}(D), 
\Vert \theta \Vert_{\mathbb R^{2}} \geq R
\} <1,
\end{equation}
and
\begin{equation} 
\label{inside-exp}
\vert \varphi (\theta) \vert \leq \exp \big( -\frac{1}{4} \langle Q\theta, \theta \rangle \big), \quad 
\Vert \theta \Vert_{\mathbb R^{2}} < R
\end{equation} hold. See e.g., \cite[P7 in Section 7]{Spitzer} and also Shiga \cite[Lemma 6.15]{Shiga} for the proof of
(\ref{inside-exp}).

By Lemma \ref{Fourier}, we have
\begin{eqnarray}
2\pi n \cdot p(n,0,y)&=&\frac{\vert {\rm det}T \vert n}{2\pi}
\int_{(\hspace{-1mm}\mbox{ }^{t}T)^{-1}(D)}
\varphi (\theta')^{n} 
e^{-\sqrt{-1} \langle y, \theta' \rangle}
d\theta' 
\nonumber \\
&=&\frac{3A(G)}{2\pi} ({\rm det}Q)^{1/2}
\int_{{\sqrt n}(\hspace{-1mm}\mbox{ }^{t}T)^{-1}(D)}
\varphi \Big (\frac{\theta}{\sqrt n} \Big)^{n}
e^{-\sqrt{-1} \langle y, \frac{\theta}{\sqrt n} \rangle }d\theta,
\label{4.04}
\end{eqnarray}
where we performed the change of variables ${\sqrt n}\theta'=\theta$ and used (\ref{Q-decompose})
and (\ref{constants})
for the second line.
We take a positive constant $r$ sufficiently small such that 
\begin{equation}
r<\min \{ C_{1}(N)K^{-\frac{2N+3}{N+3}}{\lambda}^{\frac{N}{2(N+3)}}, R \},
\label{r-condition}
\end{equation}
and divide the range of
the above integration ${\sqrt n}(\hspace{-1mm}\mbox{ }^{t}T)^{-1}(D)$
into three parts according as $\Vert \theta \Vert_{\mathbb R^{2}} \leq rn^{1/6}$; $rn^{1/6} 
<\Vert \theta \Vert_{\mathbb R^{2}} \leq R {\sqrt n}$;
$\Vert \theta \Vert_{\mathbb R^{2}} >R {\sqrt n}$. Then we can write as
%
\begin{eqnarray} 
\lefteqn{
\int_{{\sqrt n}(\hspace{-1mm}\mbox{ }^{t}T)^{-1}(D)}
\varphi \Big (\frac{\theta}{\sqrt n} \Big)^{n}
e^{-\sqrt{-1} \langle y, \frac{\theta}{\sqrt n} \rangle }d\theta}
\nonumber \\
& &
=
I(n)(y)+J_{1}(n)(y)+J_{2}(n)(y)+J_{3}(n)(y)+J_{4}(n)(y), 
\nonumber
\end{eqnarray}
where
\begin{eqnarray}
I(n)(y)&=& \sum_{j=0}^{N} I_{j}(n)(y)
:=
\sum_{j=0}^{N} \Big \{ n^{-j/2}
\int_{\mathbb R^{2}}
b_{j}(\theta)
e^{-\frac{1}{2} \langle Q \theta, \theta \rangle} 
e^{-\sqrt{-1} \langle y, \frac{\theta}{\sqrt n} \rangle }
d\theta 
\Big \},
\nonumber
\\
J_{1}(n)(y)&:=& \int_{ \Vert \theta \Vert_{\mathbb R^{2}} \leq r n^{1/6}} \Big \{ \varphi \Big (\frac{\theta}{\sqrt n} \Big)^{n}
-e^{-\frac{1}{2} \langle Q \theta, \theta \rangle} \big (\sum_{j=0}^{N} b_{j}(\theta) n^{-j/2} \big)
\Big \}
e^{-\sqrt{-1} \langle y, \frac{\theta}{\sqrt n} \rangle }
d\theta,
\nonumber \\
J_{2}(n)(y) &:=&- \int_{ \Vert \theta \Vert_{\mathbb R^{2}} > r n^{1/6}}
e^{-\frac{1}{2} \langle Q \theta, \theta \rangle} \big (\sum_{j=0}^{N} b_{j}(\theta) n^{-j/2} \big)
e^{-\sqrt{-1} \langle y, \frac{\theta}{\sqrt n} \rangle }
d\theta,
\nonumber \\
J_{3}(n)(y)&:=&\int_{r n^{1/6} <\Vert \theta \Vert_{\mathbb R^{2}} \leq R {\sqrt n}}
\varphi \Big (\frac{\theta}{\sqrt n} \Big)^{n}
e^{-\sqrt{-1} \langle y, \frac{\theta}{\sqrt n} \rangle }d\theta,
\nonumber \\
J_{4}(n)(y)&:=&\int_{\Vert \theta \Vert_{\mathbb R^{2}} 
> R{\sqrt n},~ \theta\in {\sqrt n}(\hspace{-1mm}\mbox{ }^{t}T)^{-1}(D)}
\varphi \Big (\frac{\theta}{\sqrt n} \Big)^{n}
e^{-\sqrt{-1} \langle y, \frac{\theta}{\sqrt n} \rangle }d\theta.
\nonumber
\end{eqnarray}

Our first task is to estimate the error terms $J_{1}(n)(y), J_{2}(n)(y), J_{3}(n)(y)$ and $J_{4}(n)(y)$. 
By using (\ref{outside-est})
and (\ref{inside-exp}),
we have
\begin{equation} \sup_{y\in V} \vert J_{4}(n)(y) \vert \leq \int_{{\sqrt n}(\hspace{-1mm}\mbox{ }^{t}T)^{-1}(D)} \eta^{n} d\theta
=(2\pi)^{2}\vert {\rm det}T \vert \cdot n\eta^{n},
\nonumber 
\end{equation}
and 
\begin{equation}
\sup_{y\in V} \vert J_{3}(n)(y) \vert \leq 
\int_{\Vert \theta \Vert_{\mathbb R^{2}} >rn^{1/6}} e^{-\frac{1}{4} \langle \theta, Q\theta \rangle } d\theta
\leq  \int_{\Vert \theta \Vert_{\mathbb R^{2}} >rn^{1/6}} e^{-\frac{\lambda}{4} \Vert \theta \Vert_{\mathbb R^{2}}^{2} } d\theta
\nonumber
\end{equation} 
Thus both $J_{3}(n)(y)$ and $J_{4}(n)(y)$ converge to zero as $n \to \infty$ exponentially fast uniformly for all $y\in V$. 
Similarly, we have
\begin{equation*}
\sup_{y\in V} \vert J_{2}(n)(y) \vert \leq 
e^{-\frac{\lambda}{4}n^{1/4}} \sum_{j=0}^{N} n^{-j/2} \Big(
\int_{\mathbb R^{2}} e^{-\frac{\lambda}{4} \Vert \theta \Vert^{2}_{\mathbb R^{2}}}
\vert b_{j}(\theta) \vert d\theta \Big),
\end{equation*}
and since each $b_{j}(\theta)$ has polynomial growth, we also see that
$J_{2}(n)(y)$ converges to zero as $n \to \infty$ exponentially fast uniformly for all $y\in V$.

By (\ref{Mm-est}), we have 
$$\sum_{e\in E_{0}}p(e) \langle Q^{-1}e,e \rangle^{N/2} \leq K^{N}\lambda^{-N/2}, \quad N\in \mathbb N \cup \{0\}.$$
Furthermore recalling (\ref{r-condition}), we observe that $\Vert \theta \Vert_{\mathbb R^{2}} \leq rn^{1/6}$ implies
\begin{eqnarray*}
\langle Q\theta, \theta \rangle^{1/2}
 & \leq &  K\Vert \theta \Vert_{\mathbb R^{2}}
 \nonumber \\
&\leq &
C_{1}(N)K^{-\frac{N}{N+3}}{\lambda}^{\frac{N}{2(N+3)}} n^{\frac{1}{6}}
\nonumber \\
&\leq &C_{1}(N) \Big ( \sum_{e\in E_{0}}p(e) \langle Q^{-1}e,e \rangle^{N/2} \Big )^{-\frac{1}{N+3}} n^{\frac{N+1}{2(N+3)}}.
\end{eqnarray*}
Thus we may apply Proposition \ref{phi-Taylor}, and we obtain
\begin{eqnarray}
\sup_{y\in V} \vert J_{1}(n)(y) \vert &\leq &
C_{2}(N)n^{-\frac{N+1}{2}}
\int_{\Vert \theta \Vert_{\mathbb R^{2}} \leq r n^{1/6}} 
e^{-\frac{1}{4} \langle Q\theta, \theta \rangle }
\big( \langle Q\theta, \theta \rangle^{N+3}+\langle Q\theta, \theta \rangle^{3(N+1)} \big) d\theta
\nonumber \\
&\leq &
O_{N}(n^{-\frac{N+1}{2}}). 
\nonumber
\end{eqnarray}

Now, we calculate the principal terms 
\begin{eqnarray}
I_{j}(n)(y)&:=&n^{-j/2} \int_{\mathbb R^{2}}
b_{j}(\theta)
e^{-\frac{1}{2} \langle Q \theta, \theta \rangle} 
e^{-\sqrt{-1} \langle y, \frac{\theta}{\sqrt n} \rangle }
d\theta,
\nonumber 
\quad y\in V,~j=0,1,\ldots, N. 
\end{eqnarray}
It follows directly from (\ref{Gauss-Fourier}) that
\begin{equation} I_{0}(n)(y)=2\pi ({\rm det}Q)^{-1/2} \exp \big(-\frac{1}{2n} \langle Q^{-1}y,y \rangle \big).
\nonumber
\end{equation}
Applying Proposition \ref{3.3}, we obtain
\begin{eqnarray}
I_{1}(n)(y)&=&({\sqrt -1})^{3} n^{-1/2}\int_{\mathbb R^{2}} \Big( \frac{M_{3}(\theta)}{6} \Big)
e^{-\frac{1}{2} \langle Q \theta, \theta \rangle} 
e^{-{\sqrt{-1}} \langle y, \frac{\theta}{\sqrt n} \rangle} d\theta
\nonumber \\
&=& 
2\pi n^{-1/2} ({\rm det}Q)^{-1/2} 
\exp \big( -\frac{1}{2n} \langle Q^{-1} y, y \rangle \big) 
P_{1}\big(\frac{y}{\sqrt n} \big), 
\nonumber
\end{eqnarray}
where 
\begin{equation*}
P_{1}(y):=\frac{\kappa ({\sqrt -1})^{6}}{6}
\big(  G(1,1,1)(y)
-G(2,2,2)(y)+G(3,3,3)(y) \big).
\end{equation*}
Besides, it follows from Remark \ref{3.4} that
\begin{equation}
G(i,i,i)\big(\frac{y}{\sqrt n} \big)=-n^{-3/2}\langle Q^{-1}{\bf e}_{i}, y \rangle^{3}
+3n^{-1/2} \langle Q^{-1}{\bf e}_{i}, {\bf e}_{i} \rangle \langle Q^{-1}{\bf e}_{i},y \rangle,
\quad i=1,2,3.
\label{G-iii}
\nonumber
\end{equation}
Thus 
the explicit form of $P_{1}\big(\frac{y}{\sqrt n} \big)$ is given by
\begin{eqnarray}
P_{1}\big(\frac{y}{\sqrt n} \big)
&=& \frac{\kappa}{6}n^{-3/2}
\Big \{ \langle Q^{-1}{\bf e}_{1}, y \rangle^{3}-\langle Q^{-1}{\bf e}_{2}, y \rangle^{3}+\langle Q^{-1}{\bf e}_{3}, y \rangle^{3}
\Big \}
\nonumber \\
&\mbox{ }&
+\frac{\kappa}{2}n^{-1/2} 
\Big \{ \langle Q^{-1}{\bf e}_{1}, {\bf e}_{2} \rangle \langle Q^{-1}{\bf e}_{3}, y \rangle
+\langle Q^{-1}{\bf e}_{2}, {\bf e}_{3} \rangle \langle Q^{-1}{\bf e}_{1}, y \rangle
\nonumber \\
&\mbox{ }&
\qquad
+
\langle Q^{-1}{\bf e}_{3}, {\bf e}_{1} \rangle \langle Q^{-1}{\bf e}_{2}, y \rangle
 \Big \}.
\label{P-1-1}
\end{eqnarray}

It follows from Remark \ref{R-3-9} that
\begin{eqnarray*}
I_{2}(n)(y)
&=&
n^{-2/2} \int_{\mathbb R^{2}}
({\sqrt -1})^{4} \Big( \frac{M_{4}(\theta)}{24}-\frac{M_{2}(\theta)^{2}}{8} \Big)
e^{-\frac{1}{2} \langle Q \theta, \theta \rangle} 
e^{-{\sqrt{-1}} \langle y, \frac{\theta}{\sqrt n} \rangle} d\theta
\nonumber \\
&\mbox{ }&
+
n^{-2/2} \int_{\mathbb R^{2}}
({\sqrt -1})^{6} \Big( \frac{M_{3}(\theta)^{2}}{72} \Big)
e^{-\frac{1}{2} \langle Q \theta, \theta \rangle} 
e^{-{\sqrt{-1}} \langle y, \frac{\theta}{\sqrt n} \rangle} d\theta.
\end{eqnarray*}
Then by repeating the same argument as above, we have
\begin{equation*}
I_{2}(n)(y)
=2\pi n^{-2/2} ({\rm det}Q)^{-1/2} 
\exp \big( -\frac{1}{2n} \langle Q^{-1} y, y \rangle \big) 
\Big( P_{2}^{(1)}(n)\big(\frac{y}{\sqrt n} \big)+P_{2}^{(2)}(n)\big(\frac{y}{\sqrt n} \big)
\Big),
\end{equation*}
where
\begin{eqnarray*}
%
P_{2}^{(1)}(n)(y)
&=&
\frac{(\sqrt -1)^{8}}{24} 
\big \{ {\widehat \alpha} G(1,1,1,1)(y)
+{\widehat \beta} G(2,2,2,2)(y)
\nonumber \\
&\mbox{ }&
\quad
+{\widehat \gamma} G(3,3,3,3)(y)
\big \}
\nonumber \\
&\mbox{ }&
-
\frac{(\sqrt -1)^{8}}{8} 
\big \{ {\widehat \alpha}^{2} G(1,1,1,1)(y)
+{\widehat \beta}^{2} G(2,2,2,2)(y)
\nonumber \\
&\mbox{ }&
\quad
+{\widehat \gamma}^{2}  G(3,3,3,3)(y)
+2{\widehat \alpha} {\widehat \beta} G(1,1,2,2)(y)
\nonumber \\
&\mbox{ }&
\quad
+2{\widehat \beta} {\widehat \gamma} G(2,2,3,3)(y)
+2{\widehat \gamma} {\widehat \alpha} G(1,1,3,3)(y)
\big \},
\\
P_{2}^{(2)}(n)(y)
&=&
\frac{{\kappa}^{2}(\sqrt -1)^{12}}{72} \big \{ \sum_{i=1}^{3} G(i,i,i,i,i,i)(y)
-2G(1,1,1,2,2,2)(y)
\nonumber \\
&\mbox{ }& \quad 
+2G(1,1,1,3,3,3)(y)
-2G(2,2,2,3,3,3)(y)
\big \}.
\end{eqnarray*}
Besides, it follows from (\ref{Recursive}) that
\begin{eqnarray}
G(i,i,i,j,j,j)_{(0)}(y)&=&-9 \langle Q^{-1}{\bf e}_{i}, {\bf e}_{i} \rangle
\langle Q^{-1}{\bf e}_{i}, {\bf e}_{j} \rangle
\langle Q^{-1}{\bf e}_{j}, {\bf e}_{j} \rangle
-6\langle Q^{-1}{\bf e}_{i}, {\bf e}_{j} \rangle^{3},
\nonumber
\\
G(i,i,j,j)_{(0)}(y)&=&\langle Q^{-1}{\bf e}_{i}, {\bf e}_{i} \rangle
\langle Q^{-1}{\bf e}_{j}, {\bf e}_{j} \rangle
+2\langle Q^{-1}{\bf e}_{i}, {\bf e}_{j} \rangle^{2}, \quad i,j=1,2,3.
\nonumber
\end{eqnarray}
%
By combining these identities with Lemma \ref{naiseki-keisan}, it holds
\begin{eqnarray}
& & 3 \big ( {\widehat \alpha}^{2} \langle Q^{-1}{\bf e}_{1}, {\bf e}_{1} \rangle^{2}
+{\widehat \beta}^{2} \langle Q^{-1}{\bf e}_{2}, {\bf e}_{2} \rangle^{2}
+{\widehat \gamma}^{2} \langle Q^{-1}{\bf e}_{3}, {\bf e}_{3} \rangle^{2}
\big )
\nonumber \\
& & \quad
+2 {\widehat \alpha}{\widehat \beta} \big( \langle Q^{-1}{\bf e}_{1}, {\bf e}_{1} \rangle
\langle Q^{-1}{\bf e}_{2}, {\bf e}_{2} \rangle
+2\langle Q^{-1}{\bf e}_{1}, {\bf e}_{2} \rangle^{2} \big)
\nonumber \\
& & \quad
+2 {\widehat \gamma}{\widehat \alpha} \big( \langle Q^{-1}{\bf e}_{1}, {\bf e}_{1} \rangle
\langle Q^{-1}{\bf e}_{3}, {\bf e}_{3} \rangle
+2\langle Q^{-1}{\bf e}_{1}, {\bf e}_{3} \rangle^{2} \big)
\nonumber \\
& & \quad
+2 {\widehat \beta}{\widehat \gamma} \big( \langle Q^{-1}{\bf e}_{2}, {\bf e}_{2} \rangle
\langle Q^{-1}{\bf e}_{3}, {\bf e}_{3} \rangle
+2\langle Q^{-1}{\bf e}_{2}, {\bf e}_{3} \rangle^{2} \big)
=8.
\label{P-2-1}
\end{eqnarray}
Thus the constant term of $P_{2}^{(1)}$ and $P_{2}^{(2)}$ are obtained by
\begin{equation} 
-1+
\frac{1}{8} \Big ( {\widehat \alpha} \langle Q^{-1}{\bf e}_{1}, {\bf e}_{1} \rangle^{2}
+{\widehat \beta} \langle Q^{-1}{\bf e}_{2}, {\bf e}_{2} \rangle^{2}
+
{\widehat \gamma} \langle Q^{-1}{\bf e}_{3},{\bf e}_{3} \rangle^{2}
\Big ) 
\label{P-2-2}
\end{equation}
and
\begin{eqnarray}
& &
-\frac{5}{24} \Big ( \langle Q^{-1}{\bf e}_{1}, {\bf e}_{1} \rangle^{3}+
\langle Q^{-1} {\bf e}_{2}, {\bf e}_{2} \rangle^{3}
+
\langle Q^{-1}{\bf e}_{3}, {\bf e}_{3} \rangle^{3}
\Big )
\nonumber \\
& &
+\frac{1}{6}
\Big ( \langle Q^{-1}{\bf e}_{1}, {\bf e}_{2} \rangle^{3}+
\langle Q^{-1}{\bf e}_{2}, {\bf e}_{3} \rangle^{3}
-
\langle Q^{-1}{\bf e}_{3}, {\bf e}_{1} \rangle^{3}
\Big )
\nonumber \\
& &
+\frac{1}{4}
\Big (
 \langle Q^{-1}{\bf e}_{1}, {\bf e}_{1} \rangle
\langle Q^{-1}{\bf e}_{1}, {\bf e}_{2} \rangle
\langle Q^{-1}{\bf e}_{2}, {\bf e}_{2} \rangle
+
\langle Q^{-1}{\bf e}_{2}, {\bf e}_{2} \rangle
\langle Q^{-1}{\bf e}_{2}, {\bf e}_{3} \rangle
\langle Q^{-1}{\bf e}_{3}, {\bf e}_{3} \rangle
\nonumber \\
& & \qquad
-\langle Q^{-1}{\bf e}_{3}, {\bf e}_{3} \rangle
\langle Q^{-1}{\bf e}_{3}, {\bf e}_{1} \rangle
\langle Q^{-1}{\bf e}_{1}, {\bf e}_{1} \rangle
\Big ),
\label{P-2-3}
\end{eqnarray}
respectively.
Summarizing (\ref{P-1-1}), (\ref{P-2-1}), (\ref{P-2-2}) and (\ref{P-2-3}),
we conclude that the explicit form of the coefficient of the leading term is given by
(\ref{explicit}).

For general $j=1,2,\ldots, N$, by virtue of
(\ref{BJ}), 
we may write $I_{j}(n)(y)$ as
\begin{eqnarray}
I_{j}(n)(y)&=&\sum_{k=1}^{j}I_{j}^{(k)}(n)(y)
\nonumber \\
&:=&
\sum_{k=1}^{j} \Big \{ n^{-j/2}({\sqrt -1})^{j+2k} \int_{\mathbb R^{2}} b_{j}^{(j+2k)}(\theta)
e^{-\frac{1}{2} \langle Q \theta, \theta \rangle} 
e^{-\sqrt{-1} \langle y, \frac{\theta}{\sqrt n} \rangle }d\theta \Big \}.
\nonumber 
\end{eqnarray}
Applying Proposition \ref{3.3} again, we obtain
\begin{equation}
I^{(k)}_{j}(n)(y)=2\pi ({\sqrt -1})^{j+2k} n^{-j/2} ({\rm det}Q)^{-1/2} 
\exp \big( -\frac{1}{2n} \langle Q^{-1} y, y \rangle \big) 
P_{j}^{(j+2k)}\big(\frac{y}{\sqrt n} \big),
\label{4.18}
\end{equation}
where
\begin{equation*}
P_{j}^{(j+2k)}(y)
:=
\begin{cases}
\displaystyle{\sum_{l=0}^{\frac{j-1}{2}+k}({\sqrt -1})^{j+2k-2l}P_{j,j+2k-2l}^{(j+2k)}
(y)}
& \text{ (if $j$ is odd)},
\\ 
\displaystyle{\sum_{l=0}^{\frac{j}{2}+k}({\sqrt -1})^{j+2k-2l}P_{j,j+2k-2l}^{(j+2k)}
(y)}
& \text{ (if $j$ is even) },
\end{cases}
\end{equation*}
and each 
$P_{j,j+2k-2l}^{(j+2k)}(y)$ is a homogeneous polynomial of degree $(j+2k-2l)$ in the
variables $y_{1},y_{2}$. Noting 
$(j+2k-2l) \leq 3j$ and $({\sqrt -1})^{j+2k}({\sqrt -1})^{j+2k-2l}\in \mathbb R$ 
for any $j,k,l$, we see that
\begin{equation}
P_{j}(y):=\sum_{k=1}^{j}  ({\sqrt -1})^{j+2k} 
P_{j}^{(j+2k)}(y), \quad j=1,2,\ldots, N
\label{PJ-def}
\end{equation}
is a real valued polynomial of degree at most $3j$ in the variables $y_{1},y_{2}$ and it is an odd or even function
depending on whether $j$ is odd or even.

Then by (\ref{4.18}) and (\ref{PJ-def}), we obtain
\begin{equation}
I_{j}(n)(y)=2\pi n^{-j/2} ({\rm det}Q)^{-1/2} 
\exp \big( -\frac{1}{2n} \langle Q^{-1} y, y \rangle \big) 
P_{j} \big(\frac{y}{\sqrt n} \big).
\label{IJNY}
\end{equation}
Here we mention that the term $n^{-j/2} P_{j} \big(\frac{y}{\sqrt n} \big)$ 
on the right-hand side of (\ref{IJNY})
is regarded 
as a polynomial of degree at most $2j$ in the variable $n^{-1}$.

Plugging the above all arguments into (\ref{4.04}), 
we finally obtain the desired asymptotic expansion formula (\ref{asymptotic-general}).
This completes the proof of (1).
\vspace{2mm} 

Next, we prove (2). 
For simplicity, we only give the proof in the case $k=l=0$. (In other cases, the proof
goes through in a very similar way with a slight modification.) 
Let $p_{(3)}:=p*p*p$ and we denote
by $\varphi_{(3)}$ and $Q_{(3)}$
the characteristic function and the covariance matrix associated with the $1$-step probability distribution $p_{(3)}$,
respectively. 
We define $\chi_{q}^{(3)}(\theta)$, $q\in \mathbb N$ in the same way as in
in (\ref{cum}) with $\varphi$ replaced by $\varphi_{(3)}$.
Then we easily see $Q_{(3)}=3Q$ and $\chi_{q}^{(3)}(\theta)=3\chi_{q}(\theta)$, $q\in \mathbb N$.

Let ${\widetilde G}=({\widetilde V}, {\widetilde E})$ be an enlarged triangular lattice 
defined by ${\widetilde V}:=V_{0}$ and 
${\widetilde E}:=\left\{ (x,y) \in V_{0}\times V_{0} \vert~ 
x-y \in \left\{ \pm \widetilde{\bf e}_1, \pm \widetilde{\bf e}_2, 
\pm (\widetilde{\bf e}_2-\widetilde{\bf e}_1) \right\} \right\}$ (see Figure \ref{periodic}).
Note that $$x_{1} \widehat{\bf e}_{1}+x_{2}{\widehat{\bf e}}_{2}=
y_{1}{\bf e}_{1}+y_{2}{\bf e}_{2}=z_{1}{\widetilde{\bf e}}_{1}+z_{2}{\widetilde {\bf e}}_{2}$$
implies 
\begin{equation}
\left(
       \begin{array}{cc}
       x_{1}\\
       x_{2}
       \end{array}
     \right)
     = T
     \left(
       \begin{array}{cc}
       y_{1}\\
       y_{2}
       \end{array}
     \right),
     \quad
\left(
       \begin{array}{cc}
       y_{1}\\
       y_{2}
       \end{array}
     \right)
=S
\left(
       \begin{array}{cc}
       z_{1}\\
       z_{2}
       \end{array}
     \right)
:=
\left(
       \begin{array}{cc}
       2 & 1\\
       -1& 1
       \end{array}
     \right)
\left(
       \begin{array}{cc}
       z_{1}\\
       z_{2}
       \end{array}
     \right).
\label{matrix-S}
\end{equation}

Now, we consider a
random walk on ${\widetilde G}$ whose $1$-step transition probability distribution is given by $p_{(3)}$. 
We denote by $p_{(3)}(m,x,y)$ the $m$-step transition probability of the random walk. Noting
$p_{(3)}(m,x,y):=p(3m,x,y)$, $x,y\in {\widetilde V}$
and combining Lemma \ref{Fourier} with (\ref{matrix-S}), we have
\begin{eqnarray}
2\pi (3m) \cdot p(3m,0,y)
&=& \frac{3m}{2\pi} \vert {\rm det}(TS) \vert
\int_{\{\hspace{-1mm}\mbox{ }^{t}(TS)\}^{-1}(D)}
\varphi_{(3)} (\theta)^{m} 
e^{-\sqrt{-1} \langle y, \theta \rangle}
d\theta
\nonumber \\
&=&\frac{9A(G)}{2\pi} ({\rm det}Q)^{1/2} \vert {\rm det}S \vert
\nonumber \\
&\mbox{  }&
\times 
\int_{{\sqrt m} \{\hspace{-1mm}\mbox{ }^{t}(TS)\}^{-1}(D) }
\varphi_{(3)} \Big (\frac{\theta}{\sqrt m} \Big)^{m}
e^{-\sqrt{-1} \langle y, \frac{\theta}{\sqrt m} \rangle }d\theta,
\label{4-19}
\end{eqnarray}
Here, we observe the random walk is aperiodic
because of $p_{(3)}(1,0,0)=2/9>0$.
By virtue of Spitzer \cite[P8 in Section 7]{Spitzer}, we see 
$\vert \varphi_{(3)}(\theta) \vert =1, \theta \in \{\hspace{-1mm}\mbox{ }^{t}(TS)\}^{-1}(D) $
implies $\theta=0$. Thus we may follow the proof of (1),  
and by (\ref{4-19}), we obtain
\begin{eqnarray}
2\pi (3m) \cdot p(3m,0,y)
&=& \frac{9A(G)}{2\pi} ({\rm det}Q)^{1/2} \vert {\rm det}S \vert
\sum_{j=0}^{N} {\widetilde I}_{j}(m)(y)
+O_{N}(m^{-\frac{N+1}{2}})
\label{4-20}
\end{eqnarray}
as $m\to \infty$, where
$$ {\widetilde I}_{j}(m)(y):=m^{-j/2} \int_{\mathbb R^{2}}
{\widetilde b}_{j}(\theta)
e^{-\frac{1}{2} \langle Q_{(3)} \theta, \theta \rangle} 
e^{-\sqrt{-1} \langle y, \frac{\theta}{\sqrt m} \rangle}
d\theta, \quad y\in {\widetilde V}, j=0,1,\ldots, N$$
and 
${\widetilde b}_{j}(\theta)$ 
is defined in the same way as (\ref{BJ}) with $\chi_{q}(\theta)$
replaced by $\chi_{q}^{(3)}(\theta)$.

Furthermore, repeating the calculation of $I_{j}(n)(y)$ 
and recalling $Q_{(3)}=3Q$, we obtain
\begin{equation}
{\widetilde I}_{j}(m)(y)=2\pi (3^{2}{\rm det}Q)^{-1/2} \exp \big \{ -\frac{1}{2(3m)} \langle Q^{-1}y,y \rangle \big \}
{P}_{j} \big(\frac{y}{\sqrt{3m}} \big).
\label{tilde-P}
\end{equation}
We finally obtain the desired asymptotic expansion formula
(\ref{asymptotic-3})
by combining (\ref{tilde-P}) with (\ref{4-20}) and noting $\vert {\rm det}S \vert=3$.
This completes the proof.
\qed
\subsection{Proof of Theorem \ref{FTCLT}}
First, let us suppose that $0\leq \kappa < 1/3$.
Recall that the sequence $\{ \delta_n \}_{n=1}^{\infty}$ satisfies 
$n\delta_n^2\rightarrow 3t$ as $n\to \infty$.
Since the area of the fundamental domain in 
$\mathbb{R}^2$ with respect to $\mathbb{Z}^2$-action of
$\delta_n G$ is $\delta_n^2 A(G)$, 
\begin{align*}
H_t f(x) &=\lim_{n\rightarrow \infty}
\frac{1}{2\pi t } 
\sum_{y\in V} e^{-\frac{1}{2t} \| \delta_n x_n- \delta_n y\|^2_{\mathbb{R}^2}}
 f(\delta_n y) \delta_n^2 A(G) \\
&= \lim_{n\rightarrow \infty}
\frac{3A(G)}{2\pi  n  }
\sum_{y \in V} e^{-\frac{3}{2n} \| x_n -y \|^2_{\mathbb{R}^2}
 } f(\delta_n y).
\end{align*}
Theorem \ref{Main} implies that, for any $\varepsilon>0$, there exists $N \in \mathbb{N}$ such that 
for all $n \geq N$, $|n\delta_n^2 -3t| <t$ and 
\begin{equation*}
\left|
\sum_{y \in V}
\left( 
\frac{3A(G)}{2\pi  n }
 e^{-\frac{3}{2n} \| x_n -y \|^2_{\mathbb{R}^2}  } 
 -p(n,x_n, y) \right)f(\delta_n y) \right| \leq 
\frac{3A(G)\varepsilon}{2\pi  n} \sum_{y \in V} |f(\delta_n y)|.
\end{equation*}
Since $f$ has a compact support, by the definition of $A(G)$ in (\ref{constants}) and
the property of $\delta_n$, 
\begin{align*}
\frac{3A(G)\varepsilon}{2\pi  n} \sum_{y \in V} |f(\delta_n y)|
=&\frac{3\varepsilon}{2\pi n \delta_n^2} \sum_{ y \in V} 
|f(\delta_n y) | \delta_n^2 A(G) \\
\leq &\frac{3\varepsilon}{2\pi \left( 3t -|n \delta_n^2-3t| \right) } \sum_{ y \in V} 
|f(\delta_n y) | \delta_n^2 A(G) \\
=&\frac{3\varepsilon}{ 4 \pi t }
\left(  \sum_{ y \in V} 
|f(\delta_n y) | \delta_n^2 A(G)
-\int_{\mathbb{R}^2} | f(z) |dz \right) 
+\frac{3\varepsilon}{ 4 \pi t } 
\int_{\mathbb{R}^2} |f(z)| dz 
\\
\leq & 
C\varepsilon
\end{align*}
for some $C>0$. Then we conclude Theorem \ref{FTCLT} in the case $0\leq \kappa < 1/3$.

Next, let us consider the case $\kappa=1/3$. For $n \in \mathbb{N}$, let 
$m_n \in \mathbb{N}$ 
be the quotient of $n$ divided by $3$ and $k_n \in \{ 0, 1, 2 \}$
 be its reminder. Namely, $m_n$ and $k_n$ satisfy $n=3m_n+k_n$. 
For the sequence $\{ x_n \}_{n=1}^{\infty} \subset V$ satisfying 
$\delta_n x_n \rightarrow x \in \mathbb{R}^2$ as $n\to \infty$, 
let $\{l_n \}_{n=1}^{\infty}$ be a sequence of 
$0, 1, 2$ satisfying $x_n  \in V_{l_{n}}$, a subset of $V$ given in (\ref{submodule}). 
By the periodicity of the random walk, we have
\begin{align*}
L^n (f \circ \delta_n) (x_n) = &
\sum_{ y \in V_{k_n+l_n}} p(n, x_n, y) f(\delta_n x_n)  \\
=&\sum_{y \in V_{k_n+l_n}} \left( 
p(n, x_n, y)-\frac{9A(G)}{2\pi n} e^{ -\frac{3}{2n} \| x_n-y \|_{\mathbb{R}^2}^2} \right)
f(\delta_n y ) \\
&+\frac{3}{2\pi n}
\sum_{ y \in V_{k_n+l_n}} e^{ -\frac{3}{2n} \| x_n-y \|_{\mathbb{R}^2}^2}
f(\delta_n y ) 3A(G),
\end{align*}
where $V_3=V_0$, $V_4=V_1$.
By using (2) of Theorem \ref{Main} and the previous argument 
for $0 \leq \kappa <1/3$, 
the first term of the right-hand side
converges to $0$ as 
$n \rightarrow 0$. Since the area of the fundamental domain in $\mathbb{R}^2$ 
with respect to $\mathbb{Z}^2$-action of $\delta_n V_{k_n+l_n}$ is $3\delta_n^2 A(G)$, 
\begin{equation*}
\lim_{n\rightarrow \infty}
\frac{3}{2\pi n\delta_n^2}
\sum_{ y \in V_{k_n+l_n}}e^{ -\frac{3}{2n} \| x_n-y \|_{\mathbb{R}^2}^{2} }
f(\delta_n y ) 3\delta_n^2 A(G)
= \frac{1}{2 \pi t} 
\int_{\mathbb{R}^2} 
\exp \Big( -\frac{ \| x-z \|_{\mathbb{R}^2}^2}{2t} \Big)
f (z) dz.
\end{equation*}
This completes the proof of Theorem \ref{FTCLT}.
\qed
\subsection{Proof of Theorem \ref{Trotter}}
First, let us prove (1), i.e., the convergence of the 
infinitesimal generator.
By using Taylor's theorem with respect to the coordinate $(x_1, x_2) \simeq 
x_1\mathbf{h}_1+x_2 \mathbf{h}_2$, we have 
\begin{eqnarray}
\frac{3}{\delta^2} \Delta_{d} (f \circ \delta) (x) 
&=& 
\frac{1}{\delta^2} \sum_{e  \in E_x }p(e)\left( f(\delta x)-f(\delta t(e)) \right)
\nonumber \\
&=&\frac{3}{\delta^2} \sum_{e \in E_x}
p(e) \Big( -\sum_{i=1,2} 
\frac{\partial f}{\partial x_i}(\delta x) \left( \delta t(e)-\delta x \right)_i
\nonumber \\
&\mbox{ }&
\quad \quad 
-\frac{1}{2}\sum_{i,j=1,2}\frac{\partial^2 f}{\partial x_i \partial x_j }
(\delta x) \left( \delta t(e) -\delta x \right)_i \left( \delta t(e) -\delta x \right)_j
+O(\delta^3) \Big)
\nonumber \\
&=&
-\frac{3}{\delta} \sum_{ i=1,2} \frac{\partial f}{\partial x_i}(\delta x)
\sum_{ e \in  E_x } p(e) e 
\nonumber \\
&\mbox{ }&
-\frac{3}{2} \sum_{i,j=1,2}
\sum_{e\in E_x} p(e) \left( t(e) -x \right)_i \left( t(e) - x \right)_j 
\frac{\partial^2 f}{\partial x_i \partial x_j }
(\delta x) 
\nonumber \\
&\mbox{ }&
+O(\delta).
\end{eqnarray}
Thanks to the zero mean condition {\bf(P2)}, the first term of the 
right-hand side vanishes. 
Since the covariance matrix $Q$ with respect to the standard realization is 
given by (\ref{standard-Q}), the second term of the right-hand side is
$-\frac{1}{2}(\Delta f)(\delta x)$. Hence 
we conclude (1).

%

The item (2) is deduced from (1) 
by applying Trotter's approximation theory \cite{Trotter} (see also \cite{Kotani}).
This completes the proof.
\qed
\section{Examples}
\begin{ex}[Kotani-Sunada \cite{KS00}, Example 2]
{\rm{
We consider the simple random walk. Namely, we consider the case
$\alpha=\alpha'=\beta=\beta'=\gamma=\gamma'=1/6$.
As we mentioned in Remark \ref{Traiangular explain}, 
the standard realization of the triangular lattice is the equilateral
triangular lattice in $\mathbb R^{2}$ each of whose edge has length ${\sqrt{2/3}}$. Then we have
\begin{eqnarray}
2\pi n \cdot p(n,x,y)
&=&
{\sqrt 3}\exp \Big( -\frac{3}{2 n} 
\| y-x \|^2_{\mathbb{R}^2} \Big)
\nonumber \\
&\mbox{ }&
\times
\Big( 1-\frac{1}{2}n^{-1}+\cdots+
\frac{1}{n^{N/2}}P_{N}\Big(\frac{y-x}{\sqrt n} \Big)
\Big) +O_{N} \big( n^{-\frac{N+1}{2}} \big), \quad N\geq 3
\nonumber 
\end{eqnarray}
as $n \to \infty$ uniformly for all $x,y \in V$.
}}
\end{ex}
\begin{ex}[\cite{Teruya}]
{\rm{
We consider the case 
$$\alpha=\beta=\gamma=(1/6)+\varepsilon,~ \alpha'=\beta'=\gamma'=(1/6)-\varepsilon
\quad (0< \varepsilon<1/6).$$ 
In this case, the standard
realization of the triangular lattice is same as the above example. However,
the corresponding random walk is non-symmetric with $\kappa =2\varepsilon$. Then we have
\begin{eqnarray}
2\pi n\cdot p(n,x,y)
&=&{\sqrt 3}\exp \Big( -\frac{3}{2 n} 
\| y-x \|^2_{\mathbb{R}^2} \Big)
\cdot \Big \{ 1+\Big (-\frac{1}{2}
-6\varepsilon^{2} \Big)n^{-1}
\nonumber \\
&\mbox{ }&
\quad                 
+\cdots+
\frac{1}{n^{N/2}}P_{N}\Big(\frac{y-x}{\sqrt n} \Big)
\Big \} +O_{N} \big( n^{-\frac{N+1}{2}} \big), \quad N\geq 3
\nonumber
\end{eqnarray}
as $n \to \infty$ uniformly for all $x=x_{1}{\bf h}_{1}+x_{2}{\bf h}_{2},
y=y_{1}{\bf h}_{1}+y_{2}{\bf h}_{2} \in V$.
}}
\end{ex}
\begin{ex}
{\rm{
We consider the case 
$$\alpha=1/4, \alpha'=1/12,~~ \beta=1/3, \beta'=1/6, ~~\gamma=1/6, \gamma'=0 \quad (\kappa=1/6).$$
In this case, 
\begin{equation*}
A(G)=\frac{2}{\sqrt{11}}, \ l=\frac{2\sqrt{2}}{\sqrt{11}}, \ 
\mathbf{h}_1= \!^t \left( \frac{2\sqrt{2}}{\sqrt{11}}, 0 \right), \
\mathbf{h}_2= \!^t \left( \frac{1}{\sqrt{22}}, \frac{1}{\sqrt{2}} \right)
\end{equation*}
(see Figure \ref{example3}).
\begin{figure}[tbph]
\begin{center}
\includegraphics[width=8cm]{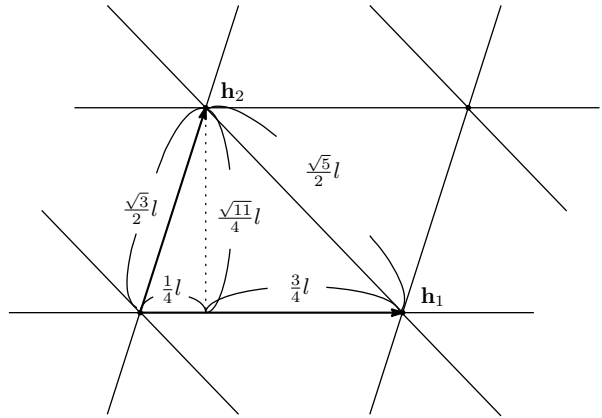} 
\end{center}
\caption{Example 5.3}
\label{example3}
\end{figure} 
Then we have
\begin{eqnarray}
2\pi n\cdot p(n,x,y)
&=&\frac{6}{ \sqrt{11} }\exp \Big( -\frac{3}{2 n} 
\| y-x \|^2_{\mathbb{R}^2} \Big)
\nonumber \\
&\mbox{ }&
\times
\Big [ 1+\Big \{ -\frac{719}{1331}
+\frac{6}{121} \Big( 2(y_{1}-x_{1})-5(y_{2}-x_{2}) \Big)
 \Big \}n^{-1} 
\nonumber \\
&\mbox{ }&
\qquad
+\cdots+
\frac{1}{{n}^{N/2}}P_{N}\Big(\frac{y-x}{\sqrt n} \Big)
\Big ] +O_{N} \big( n^{-\frac{N+1}{2}} \big), \quad N\geq 3
\nonumber
\end{eqnarray}
as $n \to \infty$ uniformly for all $x=x_{1}{\bf h}_{1}+x_{2}{\bf h}_{2},
y=y_{1}{\bf h}_{1}+y_{2}{\bf h}_{2} \in V$.
}}
\end{ex}
{\bf Acknowledgement.}
The authors would like to thank Professors Atsushi Katsuda, Motoko Kotani, 
Kazumasa Kuwada and Yukio Nagahata for 
helpful
discussions and comments.

\end{document}